\documentclass[11pt]{article}

\usepackage{amsmath,amsfonts,amsthm,amssymb,graphicx,tikz,tikz-cd,url,hyperref}
\usepackage[all,2cell,ps]{xy}

\theoremstyle{plain}
\newtheorem{thm}{Theorem}[section]

\newtheorem{lem}[thm]{Lemma}
\newtheorem{prop}[thm]{Proposition}

\newtheorem{cor}[thm]{Corollary}
\newtheorem{fact}[thm]{Fact}

\newtheorem{qtn}[thm]{Question}

\theoremstyle{definition}

\newtheorem{rem}[thm]{Remark}

\theoremstyle{remark}

\newcommand{\bbA}{\mathbb{A}}
\newcommand{\bbB}{\mathbb{B}}
\newcommand{\bbC}{\mathbb{C}}

\newcommand{\bbF}{\mathbb{F}}

\newcommand{\bbP}{\mathbb{P}}
\newcommand{\bbQ}{\mathbb{Q}}
\newcommand{\bbR}{\mathbb{R}}

\newcommand{\bbZ}{\mathbb{Z}}

\newcommand{\calA}{\mathcal{A}}

\newcommand{\calC}{\mathcal{C}}

\newcommand{\calL}{\mathcal{L}}
\newcommand{\calM}{\mathcal{M}}

\newcommand{\calO}{\mathcal{O}}

\newcommand{\calQ}{\mathcal{Q}}

\newcommand{\frakp}{\mathfrak{p}}
\newcommand{\frakq}{\mathfrak{q}}

\newcommand{\al}{\alpha}
\newcommand{\gam}{\gamma}
\newcommand{\Gam}{\Gamma}

\newcommand{\de}{\delta}
\newcommand{\Del}{\Delta}

\newcommand{\thet}{\theta}

\newcommand{\Lam}{\Lambda}
\newcommand{\sig}{\sigma}

\newcommand{\om}{\omega}
\newcommand{\Om}{\Omega}

\DeclareMathOperator{\M}{M}
\DeclareMathOperator{\SL}{SL}

\DeclareMathOperator{\GL}{GL}
\DeclareMathOperator{\PGL}{PGL}

\DeclareMathOperator{\PU}{PU}
\DeclareMathOperator{\SU}{SU}

\DeclareMathOperator{\tr}{tr}

\DeclareMathOperator{\Aut}{Aut}

\DeclareMathOperator{\Hom}{Hom}
\DeclareMathOperator{\End}{End}

\DeclareMathOperator{\Gal}{Gal}
\DeclareMathOperator{\Hol}{Hol}

\newcommand{\bs}{\backslash}

\newcommand{\conj}{\overline}
\newcommand{\wh}{\widehat}
\newcommand{\wt}{\widetilde}

\newenvironment{pf}{\begin{proof}}{\end{proof}}

\newenvironment{enum}{\begin{enumerate}}{\end{enumerate}}

\usepackage{tikz-cd}
\allowdisplaybreaks

\DeclareRobustCommand\longtwoheadrightarrow
     {\relbar\joinrel\relbar\joinrel\relbar\joinrel\twoheadrightarrow}

\makeatletter
\let\@@pmod\pmod
\DeclareRobustCommand{\pmod}{\@ifstar\@pmods\@@pmod}
\def\@pmods#1{\mkern4mu({\operator@font mod}\mkern 6mu#1)}
\makeatother

\usetikzlibrary{arrows}

\title{Algebraic fundamental groups of fake projective planes}
\author{Matthew Stover\footnote{This material is based upon work supported by Grants DMS-1906088 and DMS-2203555 from the National Science Foundation.} \\ \small{Temple University}\\ \small{\textsf{mstover@temple.edu}}}
\date{\today}

\begin{document}

\maketitle

%%%%%%%%%%%%%%%%%%%%
\begin{abstract}
Fundamental groups of fake projective planes fall into fifty distinct isomorphism classes, one for each complex conjugate pair. We prove that this is not the case for their algebraic fundamental groups: there are only forty-six isomorphism classes. We show that there are four pairs of complex conjugate pairs of fake projective planes that are $\Aut(\bbC)$-equivalent and hence have mutually isomorphic algebraic fundamental groups. All other pairs of algebraic fundamental groups are shown to be distinct through explicit finite \'etale covers. As a by-product, this provides the first examples of commensurable but nonisomorphic lattices in a rank one semisimple Lie group that have isomorphic profinite completions.
\end{abstract}
%%%%%%%%%%%%%%%%%%%%

%%%%%%%%%%%%%%%%%%%%
\section{Introduction}\label{sec:Intro}
%%%%%%%%%%%%%%%%%%%%

If $X$ is a smooth complex projective surface, then $\pi_1(X)$ will denote its topological fundamental group and $\pi_1^{\mathrm{alg}}(X)$ its algebraic fundamental group. Then $\pi_1^{\mathrm{alg}}(X)$ is the profinite completion of $\pi_1(X)$. This paper classifies the isomorphism classes of algebraic fundamental groups of fake projective planes, proving that there are exactly forty-six classes. It is known that there are fifty isomorphism classes of topological fundamental groups. Indeed, Cartwright and Steger completed the classification of fake projective planes \cite{CartwrightSteger}, proving that there are one-hundred isomorphism classes that divide into fifty distinct complex conjugate pairs. Each complex conjugate pair determines a lattice in the Lie group $\PU(2,1)$, and Mostow rigidity implies that distinct complex conjugate pairs have nonisomorphic fundamental groups. The main result of this paper is:

%%%%%%%%%%%%%%%%%%%%
\begin{thm}\label{thm:Main}
There are exactly forty-six isomorphism classes of algebraic fundamental groups of fake projective planes. Specifically, suppose that $X$ and $Y$ are fake projective planes and that $\pi_1(X) = \Gam_j, \pi_1(Y) = \Gam_k$ for $j < k$ under the numbering given in the Appendix. Then $\pi_1^{\mathrm{alg}}(X)$ is isomorphic to $\pi_1^{\mathrm{alg}}(Y)$ if and only if $\{j,k\}$ is one of the pairs
\[
\{34, 35\}\,,\,\{43, 44\}\,,\,\{47,49\}\,,\,\{48,50\}.
\]
\end{thm}
%%%%%%%%%%%%%%%%%%%%

To prove that certain fake projective planes have nonisomorphic algebraic fundamental groups, we exhibit explicit finite quotients of their topological fundamental groups that preclude the existence of an isomorphism of algebraic fundamental groups. Said in terms of \'etale covers, we produce explicit connected, \'etale, Galois covers that distinguish between the two algebraic fundamental groups. This part of the proof is the content of \S\ref{sec:Distinct}. In fact, the proof is effective. If two fake projective planes have nonisomorphic algebraic fundamental groups, we show that this can be detected by an \'etale cover of degree at most $248$; see Corollary~\ref{cor:SameList}. See Remark~\ref{rem:RefShorter} for an alternate argument found by a referee that, while ineffective and hence not sufficient to prove Corollary~\ref{cor:SameList}, cuts a much faster path to reducing the proof of Theorem~\ref{thm:Main} to proving that the four given pairs have isomorphic algebraic fundamental groups.

\medskip

The more technical and difficult part of the paper, completed in \S\ref{sec:ThePairs}, is proving that the remaining pairs have isomorphic algebraic fundamental groups. This is done following Serre's observation that $\Aut(\bbC)$ equivalent complex projective varieties have the same algebraic fundamental group \cite{Serre}. Here, we use conjugation of Shimura varieties to prove that four pairs of fake projective planes are $\Aut(\bbC)$ equivalent. Specifically, we use the close relationship between fake projective planes and Shimura varieties, which have canonical models defined over a number field, to define the action of $\Aut(\bbC)$ on fake projective planes and then apply work of Shih on conjugation of Shimura varieties \cite{Shih}; this approach has the added feature of working directly with connected components, avoiding a minor but annoying technicality that would arise using the adelic definition. This requires careful analysis of the arithmetic lattices associated with these fake projective planes. One case should also be covered by work of Borisov and Fatighenti \cite{BF} on explicit equations for fake projective planes, but this would require verification that we did not pursue (see Remark~\ref{rem:BF}). The proof of Theorem~\ref{thm:Main} therefore allows us to give the complete $\Aut(\bbC)$ classification of fake projective planes.

%%%%%%%%%%%%%%%%%%%%
\begin{cor}\label{cor:AutClassify}
There are exactly forty-six $\Aut(\bbC)$ equivalence classes of fake projective planes.
\end{cor}
%%%%%%%%%%%%%%%%%%%%

Using $\Aut(\bbC)$ conjugation to find locally symmetric varieties with nonisomorphic topological fundamental groups but isomorphic algebraic fundamental groups has some history. Examples are found in Milne and Suh \cite{MilneSuh}, Bauer, Catanese, and Grunewald \cite{BauerCataneseGrunewald}, and Rajan \cite{Rajan}. Stated in the language of discrete subgroups of Lie groups, one constructs nonisomorphic lattices in (possibly distinct) Lie groups with isomorphic profinite completions. See for example \cite{Aka, KKRS} for more on higher rank lattices with the same profinite completion.

The rank one case has recently been of particular interest, especially in light of the remarkable breakthrough by Bridson, McReynolds, Reid, and Spitler \cite{BMRS}, who found hyperbolic $3$-manifold groups $\Gam$ that are \emph{profinitely rigid}, i.e., so that if $\Lam$ is any other finitely generated, residually finite group with profinite completion $\wh{\Lam} \cong \wh{\Gam}$, then $\Lam \cong \Gam$. The first examples of nonisomorphic real rank one lattices with isomorphic profinite completions were given in \cite{StoverProfinite}, where we proved using $\Aut(\bbC)$ conjugation that, for all $n \ge 2$, there are torsion-free cocompact lattices in $\PU(n,1)$ that are not isomorphic but have isomorphic profinite completions. A distinct feature of the examples in \cite{StoverProfinite} is that the lattices are not commensurable (the higher-rank lattices in \cite{Rajan} are commensurable, but those in \cite{MilneSuh} are not). In this paper, we provide the first commensurable examples in rank one, confirming in this case the suspicion in \cite[Rem.~2.7]{MilneSuh} echoed in \cite[Rem.~2.3]{StoverProfinite}.

%%%%%%%%%%%%%%%%%%%%
\begin{cor}\label{cor:Commensurable}
There are commensurable lattices in $\PU(2,1)$ that are not isomorphic but have isomorphic profinite completions.
\end{cor}
%%%%%%%%%%%%%%%%%%%%

We note, however, that it is still an open problem whether or not there are distinct nonuniform lattices in $\PU(n,1)$ with isomorphic profinite completions. Some specific examples where certain methods used in this paper could be of use are the following:

%%%%%%%%%%%%%%%%%%%%
\begin{qtn}
Consider neat nonuniform lattices $\Gam_1, \Gam_2$ in $\PU(n,1)$ so that the associated ball quotients admit biholomorphic smooth toroidal compactifications; see for instance \cite{DiCerboStoverMultiple, DiCerboStoverCommutator, DiCerboStoverSpheres} and references therein for many (often commensurable) examples. Is it ever the case that the profinite completions $\wh{\Gam}_1$ and $\wh{\Gam}_2$ are isomorphic?
\end{qtn}
%%%%%%%%%%%%%%%%%%%%

It is noted in the Appendix to \cite{StoverVols} that the any pair of examples in \cite{DiCerboStoverClassify} with biholomorphic smooth toroidal compactifications are distinguished from one another by the number of  subgroups of index four. It follows that the associated lattices in $\PU(2,1)$ have nonisomorphic profinite completions. Thus having biholomorphic smooth toroidal compactifications does not necessarily imply an isomorphism of profinite completions.

%%%%%%%%%%%%%%%%%%%%
\begin{qtn}
Suppose that $\Gam_1, \Gam_2 < \PU(n,1)$ are neat nonuniform lattices with isomorphic profinite completions. Are the smooth toroidal compactifications of the ball quotients $\Gam_j \bs \bbB^n$ biholomorphic? Do the compactifications have isomorphic algebraic fundamental groups?
\end{qtn}
%%%%%%%%%%%%%%%%%%%%

%%%%%%%%%%%%%%%%%%%%
\begin{rem}
There are many smooth projective surfaces not birational to a fake projective plane with (topological, hence algebraic) fundamental group  isomorphic to the fundamental group of a fake projective plane. See \cite[\S 1]{StoverToledo} for a discussion of some known constructions and \cite[\S 7]{StoverToledo} for one method of producing examples. On the other hand, see \cite[\S 6]{StoverToledo} for a proof that an \emph{aspherical} smooth projective surface with fundamental group isomorphic to the fundamental group of a smooth compact ball quotient must in fact be biholomorphic to a ball quotient. In particular, an aspherical smooth projective surface with the same fundamental group as a fake projective plane must be a fake projective plane. It would be interesting to know if an aspherical smooth projective surface with $\pi_1^{\mathrm{alg}}(V)$ isomorphic to the algebraic fundamental group of a fake projective plane is necessarily a fake projective plane.
\end{rem}
%%%%%%%%%%%%%%%%%%%%

%%%%%%%%%%%%%%%%%%%%
\begin{rem}
A negative answer to \cite[Qu.~3]{StoverToledo} would provide smooth projective surfaces $V$ with $\pi_1(V)$ not isomorphic to the fundamental group of a fake projective plane but with $\pi_1^{\mathrm{alg}}(V)$ isomorphic to the algebraic fundamental group of a fake projective plane. Our later work \cite{StoverToledo2} answering the analogue of that question for arithmetic lattices in $\PU(n,1)$ of ``simple type'' does not apply to the commensurability classes that contain fake projective planes.
\end{rem}
%%%%%%%%%%%%%%%%%%%%

%%%%%%%%%%%%%%%%%%%%
\begin{rem}[Remark on use of the computer]
This paper relies on computations conducted in the computer algebra program Magma \cite{Magma}. Code is available on the author's website \cite{code}. The only functions used were those related to finite index subgroups of finitely presented groups. However, the classification of fake projective planes by Cartwright and Steger relies so heavily on computer calculation that, even if the relatively simple (if prohibitively long) computations necessary to confirm the results of this paper were done by hand, there is no reasonable sense in which this paper could be considered completely independent of the use of computer algebra software.
\end{rem}
%%%%%%%%%%%%%%%%%%%%

\begin{figure}[h]
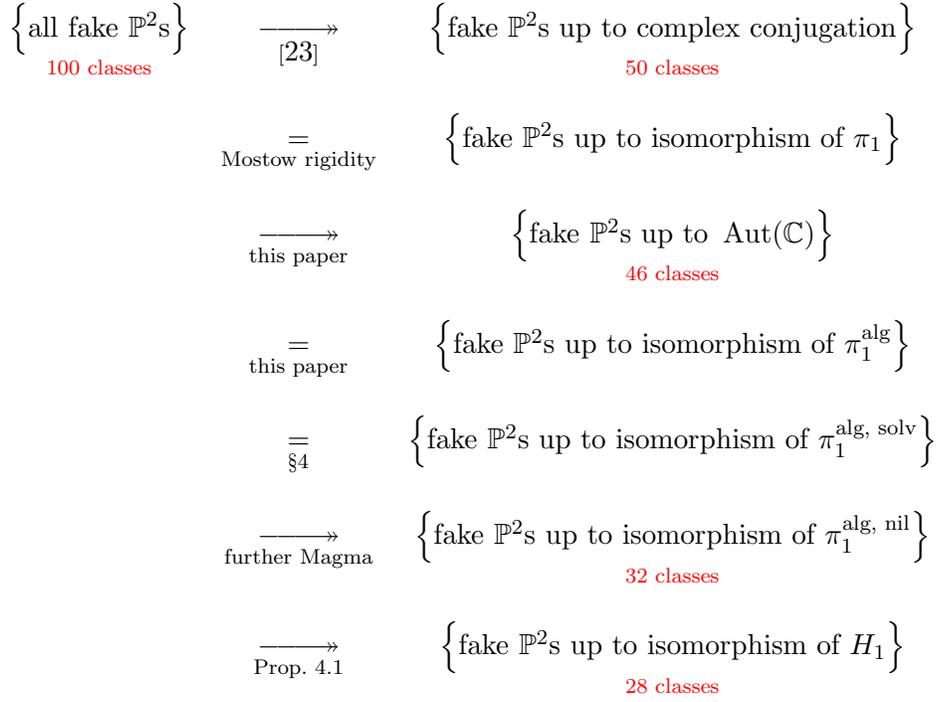

\begin{tabular}{ccc}
$\underset{\textcolor{red}{100\textrm{ classes}}}{\Big\{\textrm{all fake }\bbP^2\textrm{s}\Big\}}$ & $\underset{\cite{KulikovKharlamov}}{\longtwoheadrightarrow}$ & $\underset{\textcolor{red}{50\textrm{ classes}}}{\Big\{\textrm{fake }\bbP^2\textrm{s up to complex conjugation}\Big\}}$ \\
& & \\
& $\underset{\textrm{Mostow rigidity}}{=}$ & $\Big\{\textrm{fake }\bbP^2\textrm{s up to isomorphism of }\pi_1 \Big\}$ \\
& & \\
& $\underset{\textrm{this paper}}{\longtwoheadrightarrow}$ & $\underset{\textcolor{red}{46\textrm{ classes}}}{\Big\{\textrm{fake }\bbP^2\textrm{s up to }\Aut(\bbC)\Big\}}$ \\
& & \\
& $\underset{\textrm{this paper}}{=}$ & $\Big\{\textrm{fake }\bbP^2\textrm{s up to isomorphism of }\pi_1^{\textrm{alg}} \Big\}$ \\
& & \\
& $\underset{\textrm{\S\ref{sec:Distinct}}}{=}$ & $\Big\{\textrm{fake }\bbP^2\textrm{s up to isomorphism of }\pi_1^{\textrm{alg, solv}} \Big\}$ \\
& & \\
& $\underset{\textrm{further Magma}}{\longtwoheadrightarrow}$ & $\underset{\textcolor{red}{32\textrm{ classes}}}{\Big\{\textrm{fake }\bbP^2\textrm{s up to isomorphism of }\pi_1^{\textrm{alg, nil}} \Big\}}$ \\
& & \\
& $\underset{\textrm{Prop.\ \ref{prop:AbReduce}}}{\longtwoheadrightarrow}$ & $\underset{\textcolor{red}{28\textrm{ classes}}}{\Big\{\textrm{fake }\bbP^2\textrm{s up to isomorphism of }H_1 \Big\}}$
\end{tabular}
\caption{A hierarchy for identifying fake projective planes}\label{fig:Hi}
\end{figure}

%%%%%%%%%%%%%%%%%%%%
\noindent
\textbf{A hierarchy for fake projective plane groups}. We close the introduction by describing a hierarchy for classifying fake projective plane groups by their finite quotients, suggested by Jakob Stix. See Figure~\ref{fig:Hi}, where $\pi_1^{\mathrm{alg, solv}}$ denotes the solvable algebraic fundamental group and $\pi_1^{\mathrm{alg, nil}}$ the nilpotent. The equality between fake projective planes up to $\pi_1^{\mathrm{alg}}$ versus $\pi_1^{\mathrm{alg, solv}}$ follows from the fact that all groups used in \S\ref{sec:Distinct} to differentiate between algebraic fundamental groups are solvable.

This is not discussed elsewhere, but one can also show using routine Magma calculations that there are $32$ classes distinguished by their nilpotent quotients. Specifically, among groups $\Gam_j$ with the same abelianization, $\Gam_{41}$ is distinguished from $\Gam_4$ and $\Gam_5$ by nilpotent quotients, and similarly $\Gam_{17}$ is distinguished from $\Gam_{48}$ and $\Gam_{50}$, $\Gam_{13}$ is distinguished from $\Gam_{23}$, and $\Gam_{36}$ is distinguished from $\Gam_{37}$ and $\Gam_{45}$.

\medskip
%%%%%%%%%%%%%%%%%%%%

This paper is structured as follows. An Appendix at the end of the paper gives a numbering system to the fifty fundamental groups of fake projective planes that is referenced throughout the paper. We begin in \S\ref{sec:Action} with some basics on the action of $\Aut(\bbC)$ on smooth projective varieties. In \S\ref{sec:Profinite} we briefly present some facts about profinite completions of discrete groups that are used in \S\ref{sec:Distinct} to distinguish between all nonisomorphic algebraic fundamental groups. Then \S\ref{sec:Conj} recaps conjugation of Shimura varieties and \S\ref{sec:Buildings} briefly presents some facts about certain groups over nonarchimedean fields and their buildings. These sections are used in \S\ref{sec:ThePairs} to produce the nontrivial isomorphisms between algebraic fundamental groups that combine with the results in \S\ref{sec:Distinct} to prove Theorem~\ref{thm:Main}.

%%%%%%%%%%%%%%%%%%%
\subsubsection*{Acknowledgments} Thanks are due to Alan Reid for conversations about commensurability and isomorphisms between profinite completions, to Lev Borisov for correspondence about equations for fake projective planes, Jakob Stix for interesting correspondence, and the referees for helpful comments, particularly the shorted path to Theorem~\ref{thm:SameList} recorded in Remark~\ref{rem:RefShorter}.
%%%%%%%%%%%%%%%%%%%

%%%%%%%%%%%%%%%%%%%%
\section{The action of $\Aut(\bbC)$}\label{sec:Action}
%%%%%%%%%%%%%%%%%%%%

This section follows \cite[\S4]{BauerCataneseGrunewald}. To start, an element $\tau \in \Aut(\bbC)$ acts on an element $f(x) = \sum \al_I x^I$ of the polynomial ring $\bbC[x_0, \dots, x_N]$ by
\begin{equation*}%\label{eq:PolyAction}
\tau(f)(x) = \sum_{I = (i_0, \dots, i_N)} \tau(\al_I) x^I,
\end{equation*}
where we use the usual multi-index notation for polynomials. Given a projective variety $X \subseteq \bbP^N$ cut out by homogeneous polynomials $\{f_1(x), \dots, f_r(x)\}$, the coordinate action of $\tau \in \Aut(\bbC)$ on $\bbP^n$ maps $X$ to a new variety
\begin{equation*}%\label{eq:VarietyAction}
X^\tau = \{x \in \bbP^N\ :\ \tau(f_j)(x) = 0 \textrm{ for } 1 \le j \le r\}.
\end{equation*}
Indeed, $f_j(x) = 0$ if and only if $\tau(f_j)(\tau(x)) = 0$.

Applying the above to the graph of an element $\psi \in \Aut(X)$, we obtain an element $\psi^\tau \in \Aut(X^\tau)$. A \emph{$G$-marked variety} is then a triple $(X, G, \eta)$, where $X$ is a projective variety and $\eta : G \to \Aut(X)$ is an injective homomorphism. For each $G$-marking $(X, G, \eta)$ of $X$, we obtain the $G$-marking
\begin{equation}\label{eq:Marking}
(X, G, \eta)^\tau = (X^\tau, G, \tau \circ \eta \circ \tau^{-1})
\end{equation}
of $X^\tau$. Set $\eta^\tau = \tau \circ \eta \circ \tau^{-1}$. We record some essential properties of the $\Aut(\bbC)$ action on projective varieties.

%%%%%%%%%%%%%%%%%%%%
\begin{lem}\label{lem:AutProperties}
Let $X \subseteq \bbP^N$ be a projective variety. For $\tau \in \Aut(\bbC)$ and a marking $(X, G, \eta)$, the following hold.
\begin{enum}\itemsep-0.25em

\item If $\eta(g)$ has a fixed point for some $g \in G$, then $\eta^\tau(g)$ has a fixed point on $X^\tau$.

\item One has $(X / \eta(G))^\tau = X^\tau / \eta^\tau(G)$.

\item The variety $X$ is smooth if and only if $X^\tau$ is smooth.

\end{enum}
\end{lem}
%%%%%%%%%%%%%%%%%%%%

%%%%%%%%%%%%%%%%%%%%
\begin{pf}
The first property follows from the definition of $\eta^\tau$ and the fact that $X \to X^\tau$ is a bijection on points. The second property is \cite[Rem.~4.4(4)]{BauerCataneseGrunewald}. The third follows from the fact that smoothness is defined by the rank of the Jacobian matrix for the system of polynomials, which is invariant under the action of $\Aut(\bbC)$.
\end{pf}
%%%%%%%%%%%%%%%%%%%%

We will use the following corollary in \S\ref{sec:ThePairs}.

%%%%%%%%%%%%%%%%%%%%
\begin{cor}\label{cor:FakeQuotient}
Let $G$ be a finite group and $(Y, G, \eta)$ be a smooth projective surface such that $\eta(G)$ has no fixed points (i.e., $Y \to X$ is finite \'etale) and $X = Y / \eta(G)$ is a fake projective plane. Then $X^\tau$ is a fake projective plane.
\end{cor}
%%%%%%%%%%%%%%%%%%%%

%%%%%%%%%%%%%%%%%%%%
\begin{pf}
Lemma~\ref{lem:AutProperties} implies that $X^\tau = Y^\tau / \eta^\tau(G)$ is also a smooth projective surface. Betti numbers are invariant under $\tau$ by Serre's GAGA principle, so $X^\tau$ also has the same complex cohomology as $\bbP^2$. We must rule out the possibility that $X^\tau$ is isomorphic to $\bbP^2$. Since $X$ is a ball quotient, so is $X^\tau$ by a theorem of Kazhdan \cite{Kazhdan}. In particular, $X^\tau$ has infinite fundamental group, and therefore it is a fake projective plane.
\end{pf}
%%%%%%%%%%%%%%%%%%%%

%%%%%%%%%%%%%%%%%%%%
\section{Profinite completions of discrete groups}\label{sec:Profinite}
%%%%%%%%%%%%%%%%%%%%

Let $X$ be a smooth complex projective variety. Throughout this paper, $\pi_1(X)$ will denote the topological fundamental group, considered as a finitely generated discrete group, and $\pi_1^{\mathrm{alg}}(X)$ the algebraic fundamental group. Then $\pi_1^{\mathrm{alg}}(X)$ is isomorphic to the profinite completion of $\pi_1(X)$. As explained in the introduction, there are now many known examples of pairs $X, Y$ so that $\pi_1(X)$ is not isomorphic to $\pi_1(Y)$ but $\pi_1^{\mathrm{alg}}(X) \cong \pi_1^{\mathrm{alg}}(Y)$. Here and throughout, the profinite completion of a group $\Gam$ is denoted $\wh{\Gam}$.

This section is concerned with tools to prove that $\pi_1^{\mathrm{alg}}(X)$ is not isomorphic to $\pi_1^{\mathrm{alg}}(Y)$ through finite quotients of $\pi_1(X)$ and $\pi_1(Y)$. These tools come from the study of profinite completions of arbitrary finitely generated discrete groups. See \cite{ReidProfinite} for a thorough survey; we will only need the following two fundamental facts found there:

%%%%%%%%%%%%%%%%%%%%
\begin{lem}[Prop.~3.2 \cite{ReidProfinite}]\label{lem:ProAb}
If $\Gam$ and $\Del$ are discrete groups with $\wh{\Gam} \cong \wh{\Del}$, then $\Gam$ and $\Del$ have isomorphic abelianizations.
\end{lem}
%%%%%%%%%%%%%%%%%%%%

%%%%%%%%%%%%%%%%%%%%
\begin{thm}[Cor.~4.2 \cite{ReidProfinite}]\label{thm:QuoNum}
If $\Gam, \Del$ are finitely generated groups with $\wh{\Gam} \cong \wh{\Del}$, then one has equalities
\begin{align*}
\left|\Hom(\Gam, Q)\right| &= \left|\Hom(\Del, Q)\right| \\ %\label{eq:QuoNum} \\
\left|\mathrm{Epi}(\Gam, Q)\right| &= \left|\mathrm{Epi}(\Del, Q)\right| %\label{eq:EpiNum}
\end{align*}
for all finite groups $Q$, where $\mathrm{Epi}(\Lam, Q)$ is the set of epimorphisms of a group $\Lam$ onto $Q$.
\end{thm}
%%%%%%%%%%%%%%%%%%%%

In the language of smooth complex projective varieties, Lemma~\ref{lem:ProAb} says that if $\pi_1^{\mathrm{alg}}(X) \cong \pi_1^{\mathrm{alg}}(Y)$, then $H_1(X, \bbZ) \cong H_1(Y, \bbZ)$. Theorem~\ref{thm:QuoNum} says that if $X$ and $Y$ admit a different number of isomorphism classes of connected \'etale coverings with Galois group some fixed finite group $Q$, then they have nonisomorphic algebraic fundamental groups. Note that Theorem~\ref{thm:QuoNum} is slightly stronger than the more elementary statement that if $X$ admits a connected \'etale cover with group $Q$ and $Y$ does not, then their algebraic fundamental groups are not isomorphic. (Here we call \'etale covers isomorphic if they are associated with the same conjugacy class of subgroups of the fundamental group.)

%%%%%%%%%%%%%%%%%%%%
\section{Distinct isomorphism classes}\label{sec:Distinct}
%%%%%%%%%%%%%%%%%%%%

In this section, we prove that the distinct algebraic fundamental groups of fake projective planes are indeed distinct. Throughout this section, we use the numbering of fake projective plane fundamental groups given in the Appendix. Applying Lemma~\ref{lem:ProAb} after using computer algebra software to compute abelianizations (also see \cite{Cartwright}) immediately gives us the following start to proceedings.

%%%%%%%%%%%%%%%%%%%%
\begin{prop}\label{prop:AbReduce}
If $\Gam_j$ and $\Gam_k$ are fake projective plane fundamental groups with $\wh{\Gam}_j \cong \wh{\Gam}_k$, then $j$ and $k$ appear in the same cell of Table~\ref{tb:AbReduce}.
\begin{table}[h]
\centering
\begin{tabular}{|c|c|}
\hline
$H_1(\Gam_j, \bbZ)$ & $j$ \\
\hline
$(\bbZ / 2) \times (\bbZ / 26)$ & $4, 5, 41$ \\
\hline
$\bbZ / 14$ & $9, 11, 25, 32, 39$ \\
\hline
$(\bbZ / 2) \times (\bbZ / 14)$ & $12, 20$ \\
\hline
$(\bbZ / 2)^2 \times (\bbZ / 6)$ & $13, 23$ \\
\hline
$\bbZ / 28$ & $16, 47, 49$ \\
\hline
$(\bbZ / 2) \times (\bbZ / 12)$ & $17, 48, 50$ \\
\hline
$(\bbZ / 2) \times (\bbZ / 18)$ & $18, 26$ \\
\hline
\end{tabular}
\begin{tabular}{|c|c|}
\hline
$H_1(\Gam_j, \bbZ)$ & $j$ \\
\hline
$\bbZ / 18$ & $19, 34, 35$ \\
\hline
$\bbZ / 42$ & $22, 30$ \\
\hline
$\bbZ / 6$ & $24, 43, 44$ \\
\hline
$\bbZ / 21$ & $27, 31$ \\
\hline
$(\bbZ / 3) \times (\bbZ / 6)$ & $36, 37, 45$ \\
\hline
$\bbZ / 7$ & $38, 40$ \\
\hline
\multicolumn{1}{c}{} & \multicolumn{1}{c}{} \\
\end{tabular}
\caption{Fake projective planes with isomorphic first homology over $\bbZ$.}\label{tb:AbReduce}
\end{table}
\end{prop}
%%%%%%%%%%%%%%%%%%%%

Note that Proposition~\ref{prop:AbReduce} omits those $1 \le j \le 50$ where $\Gam_j$ is the unique fake projective plane fundamental group with abelianization isomorphic to $\Gam_j^{ab}$. Consequently, for any $j$ that does not appear in a cell of Table~\ref{tb:AbReduce}, Lemma~\ref{lem:ProAb} immediately gives us that $\wh{\Gam}_j \cong \wh{\Gam}_k$ implies $k = j$; these $\wh{\Gam}_j$ are therefore completely classified with regard to the proof of Theorem~\ref{thm:Main}.

We now focus on finite quotients that distinguish between various cells of Table~\ref{tb:AbReduce}. The following sequence of facts is proved with computer algebra software by computing the number of homomorphisms onto the appropriate finite group $G$. We use the following notation for finite groups: $S_n$ is the symmetric group on $n$ letters, $A_n$ the alternating group on $n$ letters, and $D_n$ the dihedral group of order $2 n$. In what follows, two finite quotients of a group are considered distinct if and only if they have distinct kernels.

%%%%%%%%%%%%%%%%%%%%
\begin{fact}\label{fact:SmGp6,1}
If $\Gam_j$ is a fake projective plane group admitting $S_3$ as a quotient, then
\begin{align*}
j \in \{1, &\,2, 4, 5, 11, 12, 14, 16, 20, 21, 25, 26, 30,\\
&\,32, 33, 34, 35, 36, 37, 42, 43, 44, 45, 47, 49 \}.
\end{align*}
The $S_3$ quotient is unique except in the following cases:
\[
\begin{cases}
j = 20 & 4\,\textrm{quotients} \\
j = 21 & 2\,\textrm{quotients} \\
j = 25 & 4\,\textrm{quotients} \\
j = 26 & 2\,\textrm{quotients} \\
j = 32 & 4\,\textrm{quotients}
\end{cases}
\]
\end{fact}
%%%%%%%%%%%%%%%%%%%%

%%%%%%%%%%%%%%%%%%%%
\begin{fact}\label{fact:SmGp12,3}
If $\Gam_j$ is a fake projective plane group admitting an $A_4$ quotient, then $j \in \{22, 24, 27, 29, 45\}$ and each such $\Gam_j$ admits exactly two distinct $A_4$ quotients.
\end{fact}
%%%%%%%%%%%%%%%%%%%%

%%%%%%%%%%%%%%%%%%%%
\begin{fact}\label{fact:SmGp48,29}
The only fake projective plane group with a $\GL_2(\bbF_3)$ quotient is $\Gam_5$, which admits four such quotients.
\end{fact}
%%%%%%%%%%%%%%%%%%%%

%%%%%%%%%%%%%%%%%%%%
\begin{fact}\label{fact:SmGp56,11}
If $\Gam_j$ is a fake projective plane group with the affine linear group $\mathrm{AGL}_1(\bbF_8)$ in one variable over $\bbF_8$ as a quotient, then
\[
j \in \{9, 11, 12, 14, 16, 20, 22, 25, 27, 30, 31, 32, 38, 39, 40, 47, 49\}.
\]
There are twenty-seven $\mathrm{AGL}_1(\bbF_8)$ quotients for $j \in \{20, 32, 38\}$ and three otherwise.
\end{fact}
%%%%%%%%%%%%%%%%%%%%

%%%%%%%%%%%%%%%%%%%%
\begin{fact}\label{fact:SmGp16,7}
If $\Gam_j$ is a fake projective plane admitting $D_8$ as a quotient, then $j \in \{8, 48, 50\}$ and $\Gam_j$ admits precisely four $D_8$ quotients.
\end{fact}
%%%%%%%%%%%%%%%%%%%%

%%%%%%%%%%%%%%%%%%%%
\begin{fact}\label{fact:SmGp18,1}
The only fake projective plane group admitting $D_9$ as a quotient is $\Gam_{32}$, which admits nine such quotients.
\end{fact}
%%%%%%%%%%%%%%%%%%%%

%%%%%%%%%%%%%%%%%%%%
\begin{fact}\label{fact:SmGp26,1}
The only fake projective plane group admitting $D_{13}$ as a quotient is $\Gam_{39}$, which admits six such quotients.
\end{fact}
%%%%%%%%%%%%%%%%%%%%

%%%%%%%%%%%%%%%%%%%%
\begin{fact}\label{fact:SmGp8,4}
If $\Gam_j$ is a fake projective plane group with the quaternion group $Q$ of order $8$ as a quotient, then $j \in \{3, 23, 28, 41\}$. For $j = 3, 23$, $\Gam_j$ admits twelve $Q$ quotients and for $j = 28, 41$ there are only six.
\end{fact}
%%%%%%%%%%%%%%%%%%%%

%%%%%%%%%%%%%%%%%%%%
\begin{fact}\label{fact:SmGp27,4}
If $\Gam_j$ is a fake projective plane group having the unique nonabelian group $G$ of order $27$ that is a semidirect product $(\bbZ / 9) \rtimes (\bbZ / 3)$ as a quotient, then $j \in \{37, 45\}$ and $\Gam_j$ admits six such quotients.
\end{fact}
%%%%%%%%%%%%%%%%%%%%

%%%%%%%%%%%%%%%%%%%%
\begin{fact}\label{fact:SmGp52,3}
If $\Gam_j$ is a fake projective plane group with the unique semidirect product $(\bbZ / 13) \rtimes (\bbZ / 4)$ with faithful action of $\bbZ / 4$ as a quotient, then ${j \in \{1, 47, 49\}}$ and $\Gam_j$ admits three such quotients.
\end{fact}
%%%%%%%%%%%%%%%%%%%%

We now combine the above information to reduce the proof of Theorem~\ref{thm:Main} to considering four pairs.

%%%%%%%%%%%%%%%%%%%%
\begin{thm}\label{thm:SameList}
Let $X$ and $Y$ be fake projective planes with $\pi_1(X) = \Gam_j$ and $\pi_1(Y) = \Gam_k$ for $j < k$. If the pair $\{j, k\}$ is not among
\[
\{34, 35\}\,,\,\{43, 44\}\,,\,\{47,49\}\,,\,\{48,50\},
\]
then $\pi_1^{\mathrm{alg}}(X)$ and $\pi_1^{\mathrm{alg}}(Y)$ are not isomorphic.
\end{thm}
%%%%%%%%%%%%%%%%%%%%

%%%%%%%%%%%%%%%%%%%%
\begin{pf}
We must prove $\wh{\Gam}_j$ is not isomorphic to $\wh{\Gam}_k$ with the exception of the four pairs in the statement of the theorem. First, note that $j$ and $k$ must appear in the same cell of Table~\ref{tb:AbReduce}. We now apply Theorem~\ref{thm:QuoNum} to the quotients given in Facts~\ref{fact:SmGp6,1}-\ref{fact:SmGp52,3} to eliminate the remaining pairs.

\medskip
\noindent
\textbf{Case 1}: $\{4, 5, 41\}$ (Facts~\ref{fact:SmGp6,1} and \ref{fact:SmGp48,29})
\medskip

The group $\Gam_{41}$ does not have $S_3$ as a quotient, but $\Gam_4$ and $\Gam_5$ do. Then $\Gam_4$ does not have $\GL_2(\bbF_3)$ as a quotient, but $\Gam_5$ does. Thus all three have distinct profinite completions.

\medskip
\noindent
\textbf{Case 2}: $\{9, 11, 25, 32, 39\}$ (Facts~\ref{fact:SmGp6,1}, \ref{fact:SmGp18,1}, and \ref{fact:SmGp26,1})
\medskip

The groups $\Gam_{11}$ and $\Gam_{25}$ have an $S_3$ quotient, but $\Gam_9$, $\Gam_{32}$, and $\Gam_{39}$ do not. Moreover, $\Gam_{25}$ has four $S_3$ quotients, while $\Gam_{11}$ only has one. Then $\Gam_{32}$ has a $D_9$ quotient, unlike $\Gam_9$ and $\Gam_{25}$, and $\Gam_{25}$ has a $D_{13}$ quotient, unlike $\Gam_9$ or $\Gam_{32}$. Thus all five groups have distinct profinite completions.

\medskip
\noindent
\textbf{Case 3}: $\{12, 20\}$ (Fact~\ref{fact:SmGp6,1})
\medskip

The groups $\Gam_{12}$ and $\Gam_{20}$ both have $S_3$ quotients, but $\Gam_{12}$ has only one, whereas $\Gam_{20}$ has four. Therefore the pair have distinct profinite completions.

\medskip
\noindent
\textbf{Case 4}: $\{13, 23\}$ (Fact~\ref{fact:SmGp8,4})
\medskip

The group $\Gam_{23}$ has the quaternion group of order $8$ as a quotient, but $\Gam_{13}$ does not, hence the two have nonisomorphic profinite completions.

\medskip
\noindent
\textbf{Case 5}: $\{16, 47, 49\}$ (Fact~\ref{fact:SmGp52,3})
\medskip

The unique semidirect product $(\bbZ / 13) \rtimes (\bbZ / 4)$ with faithful $\bbZ / 4$-action is a quotient of $\Gam_{47}$ and $\Gam_{49}$, but not of $\Gam_{16}$. Therefore $\wh{\Gam}_{16}$ cannot be isomorphic to $\wh{\Gam}_j$ for $j \in \{47, 49\}$.

\medskip
\noindent
\textbf{Case 6}: $\{17, 48, 50\}$ (Fact~\ref{fact:SmGp16,7})
\medskip

The groups $\Gam_{48}$ and $\Gam_{50}$ both admit $D_8$ as a quotient, but $\Gam_{17}$ does not. Therefore $\wh{\Gam}_{17}$ cannot be isomorphic to $\wh{\Gam}_j$ for $j \in \{48, 50\}$.

\medskip
\noindent
\textbf{Case 7}: $\{18, 26\}$ (Fact~\ref{fact:SmGp6,1})
\medskip

The group $\Gam_{18}$ does not admit $S_3$ as a quotient, but $\Gam_{26}$ does. Therefore the pair have distinct profinite completions.

\medskip
\noindent
\textbf{Case 8}: $\{19, 34, 35\}$ (Fact~\ref{fact:SmGp6,1})
\medskip

The groups $\Gam_{34}$ and $\Gam_{35}$ both have $S_3$ as a quotient, but $\Gam_{19}$ does not. Therefore $\wh{\Gam}_{19}$ cannot be isomorphic to $\wh{\Gam}_j$ for $j \in \{34, 35\}$.

\medskip
\noindent
\textbf{Case 9}: $\{22, 30\}$ (Fact~\ref{fact:SmGp6,1})
\medskip

The group $\Gam_{30}$ has an $S_3$ quotient, unlike $\Gam_{22}$, so the pair have nonisomorphic profinite completions.

\medskip
\noindent
\textbf{Case 10}: $\{24, 43, 44\}$ (Fact~\ref{fact:SmGp6,1})
\medskip

The groups $\Gam_{43}$ and $\Gam_{44}$ both admit $S_3$ as a quotient, but $\Gam_{24}$ does not. Therefore $\wh{\Gam}_{24}$ cannot be isomorphic to $\wh{\Gam}_j$ for $j \in \{43, 44\}$.

\medskip
\noindent
\textbf{Case 11}: $\{27, 31\}$ (Fact~\ref{fact:SmGp12,3})
\medskip

The group $\Gam_{27}$ has $A_4$ as a quotient, unlike $\Gam_{31}$, so their profinite completions are not isomorphic.

\medskip
\noindent
\textbf{Case 12}: $\{36, 37, 45\}$ (Facts~\ref{fact:SmGp12,3} and \ref{fact:SmGp27,4})
\medskip

The group $\Gam_{45}$ admits an $A_4$ quotient, but $\Gam_{36}$ and $\Gam_{37}$ do not. Also $\Gam_{37}$ has the nontrivial semidirect product $(\bbZ / 9) \rtimes (\bbZ / 3)$ as a quotient, but $\Gam_{36}$ does not. Therefore all three have distinct profinite completions.

\medskip
\noindent
\textbf{Case 13}: $\{38, 40\}$ (Fact~\ref{fact:SmGp56,11})
\medskip

Each group has $\mathrm{AGL}_1(\bbF_8)$ as a quotient, but $\Gam_{38}$ has twenty-seven distinct $\mathrm{AGL}_1(\bbF_8)$ quotients, whereas $\Gam_{40}$ has only three. It follows that their profinite completions are not isomorphic.

\medskip

We have ruled out all possible isomorphisms between algebraic fundamental groups of fake projective planes except the four pairs in the statement of the theorem. Therefore the proof of the theorem is complete.
\end{pf}
%%%%%%%%%%%%%%%%%%%%

Define two connected \'etale covers $Z \to X$ between smooth projective varieties to be isomorphic when the associated subgroups of $\pi_1(X)$ are conjugate. In particular, under this definition a fixed $X$ can admit nonisomorphic covers by a fixed $Z$. In this language, considering the orders of the finite groups in Facts~\ref{fact:SmGp6,1}-\ref{fact:SmGp52,3} and possible orders of first homology groups with $\bbZ$ coefficients, we obtain the following corollary to Theorem~\ref{thm:SameList}.

%%%%%%%%%%%%%%%%%%%%
\begin{cor}\label{cor:SameList}
Let $X$ and $Y$ be fake projective planes for which the number of isomorphism classes of \'etale covers with group $G$ are exactly the same for all groups $G$ with $|G| \le 248$. Then $X$ and $Y$ are either isomorphic, complex conjugate, or $\pi_1(X) = \Gam_j$ and $\pi_1(Y) = \Gam_k$ where (up to reordering) the pair $\{j, k\}$ is one of the four given in the statement of Theorem~\ref{thm:SameList}.
\end{cor}
%%%%%%%%%%%%%%%%%%%%

In other words, away from the four pairs of fundamental groups we will soon show have isomorphic profinite completions, one can effectively determine through \'etale covers whether or not two fake projective planes are isomorphic up to complex conjugation.

%%%%%%%%%%%%%%%%%%%%
\begin{rem}\label{rem:RefShorter}
A referee noticed the following shorter path to proving Theorem~\ref{thm:SameList}. Let $X^{\mathrm{ab}}$ be the universal \'etale abelian covering of the fake projective plane $X$. If $\pi_1^{\mathrm{alg}}(X)$ is isomorphic to $\pi_1^{\mathrm{alg}}(Y)$, then one also has $H_1(X^{\mathrm{ab}}, \bbZ) \cong H_1(Y^{\mathrm{ab}}, \bbZ)$. A short computer program combining this observation with Proposition~\ref{prop:AbReduce} reduces one to studying the pairs
\[
\{4,5\}\,,\, \{34, 35\}\,,\, \{43, 44\}\,,\, \{47, 49\}\,,\, \{48, 50\}.
\]
One then rules out the pair $\{4,5\}$ using Fact~\ref{fact:SmGp48,29}. However, this argument does not imply Corollary~\ref{cor:SameList}, and the small nonabelian \'etale covers produced in the current proof of Theorem~\ref{thm:SameList} may be of independent interest, so the longer proof was retained.
\end{rem}
%%%%%%%%%%%%%%%%%%%%

%%%%%%%%%%%%%%%%%%%%
\section{Conjugation of Shimura varieties}\label{sec:Conj}
%%%%%%%%%%%%%%%%%%%%

This section recaps work of Milne and Suh \cite[\S 1]{MilneSuh} and earlier work of Miyake \cite{Miyake} on conjugation of Shimura varieties and generally follows the notation from \cite{MilneSuh}. Let $F$ be a totally real number field and $H$ an algebraic group over $F$. The restriction of scalars $H_*$ is the $\bbQ$-algebraic group so that
\begin{equation*}%\label{eq:Restrict1}
H_*(R) = H(F \otimes_\bbQ R)
\end{equation*}
for any $\bbQ$-algebra $R$ and
\begin{equation*}%\label{eq:Restrict2}
H_*(\bbR) = \prod_{\nu : F \to \bbR} H_\nu(\bbR),
\end{equation*}
where $H_\nu$ is the real algebraic group defined by extension of scalars from $\nu(F)$ to $\bbR$. Notice that $H(F)$ is a subgroup of $H_*(\bbR)$, and we then define
\begin{equation*}%\label{eq:PlusGroup}
H(F)^+ = H(F) \cap H_*(\bbR)^+,
\end{equation*}
where $H_*(\bbR)^+$ is the connected component of $H_*(\bbR)$ containing the identity. We assume that $H$ is a simply connected semisimple $F$-algebraic group, which implies that $H(\bbR)$ is connected. If $\nu$ is a nonarchimedean place of $F$ with completion $F_\nu$, then $H_{F_\nu}$ will denote the group $H(F_\nu)$ obtained by extension of scalars from $F$ to $F_\nu$.

Let $\bbA_F^\infty$ denote the finite adeles of $F$. A \emph{congruence subgroup} of $H(F)$ is a subgroup of the form
\begin{equation*}%\label{eq:CongDef}
\wt{\Gam} = H(F) \cap K
\end{equation*}
for an open compact subgroup $K < H(\bbA_F^\infty)$, where $H(F)$ is diagonally embedded in $H(\bbA_F^\infty)$. If $H^{ad}(F)^+$ is the adjoint group of $H(F)^+$, then a subgroup $\Gam < H^{ad}(F)^+$ is a congruence subgroup if its preimage $\wt{\Gam}$ in $H(F)$ with respect to the natural central isogeny is a congruence subgroup. When we need to differentiate between different congruence subgroups, we will use $\wt{\Gam}(K)$ and $\Gam(K)$ to emphasize the open compact subgroup.

Suppose that $X$ is a bounded symmetric domain and $\Hol(X)^+$ is the connected component of the identity in the group of holomorphic automorphisms of $X$. The case of interest to us is where there is a surjective homomorphism $H(\bbR)^+ \to \Hol(X)^+$ with compact kernel. If $\Gam < H^{ad}(F)^+$ is a congruence subgroup and $V = \Gam \bs X$ the associated algebraic variety, then we say that $V$ has \emph{type} $(H, X)$.

The following is Theorem 1.3 of \cite{MilneSuh} combined with the more refined version of the theorem's conclusion provided by \cite[Rem.~1.6]{MilneSuh}.

%%%%%%%%%%%%%%%%%%%%
\begin{thm}\label{thm:MilneSuh}
Suppose that $H$ is a simply connected and semisimple $F$-algebraic group and that $V$ is a smooth algebraic variety over $\bbC$ of type $(H, X)$. For any $\tau \in \Aut(\bbC)$, the variety $V^\tau$ is of type $(H^\prime, X^\prime)$ for some semisimple $F$-algebraic group $H^\prime$ such that:
\begin{equation*}%\label{eq:MilneSuh}
\begin{cases}
H_\nu^\prime \cong H_{\tau \circ \nu} & \textrm{for all real archimedean } \nu \\
 & \\
{H^\prime}_{\hspace{-0.2em}F_\nu} \cong H_{F_\nu} & \textrm{for all nonarchimedean } \nu
\end{cases}
\end{equation*}
Moreover, suppose that $V = \wt{\Gam}(K) \bs X$ for $K \le H(\bbA_F^\infty)$ with $\wt{\Gam}(K)$ torsion-free. The isomorphisms ${{H^\prime}_{\hspace{-0.2em}F_\nu} \cong H_{F_\nu}}$ for all nonarchimedean $\nu$ induce an isomorphism ${H^\prime(\bbA_F^\infty) \cong H(\bbA_F^\infty)}$ so that
\begin{equation*}%\label{eq:MilneSuh2}
V^\tau = \wt{\Gam}^\prime(K^\prime) \bs X^\prime,
\end{equation*}
where $K^\prime$ is the image of $K$ in $H^\prime(\bbA_F^\infty)$ under the isomorphism.
\end{thm}
%%%%%%%%%%%%%%%%%%%%

In fact, we will need even more concrete information about the induced isomorphism $H(\bbA_F^\infty) \to H^\prime(\bbA_F^\infty)$, so we sketch the ideas behind the proof of Theorem~\ref{thm:MilneSuh}. We follow Miyake \cite[\S5]{Miyake}, which specifically considers the case of interest for us.

\medskip

Let $F$ be a totally real number field with $[F : \bbQ] = d$, $\calO_F$ be the ring of integers of $F$, and $E / F$ be a totally imaginary quadratic extension. Fix a cubic division algebra $D$ with center $E$ and involution of second kind, that is, an anti-involution $\sig$ so that the restriction of $\sig$ to $E$ is the nontrivial element of $\Gal(E/F)$. A nonzero element $h \in D^*$ for which $\sig(h) = h$ defines a nondegenerate hermitian form on $D$ by
\begin{equation*}%\label{eq:GenHermitian}
h(x, y) = \sig(x) h y.
\end{equation*}
One then modifies $h$ to obtain a symplectic form $\psi$ on $D$ (as an $F$-vector space) in a standard way described in \cite[Ex.~8.5]{Milne}.

Set $n = 9 d$. Considering $D_\bbR = D \otimes_\bbQ \bbR$ as a complex vector space, the action of $D$ on itself defines a representation $\Phi : D \to \M_n(\bbC)$ so that $\Phi \oplus \conj{\Phi}$ is a rational representation, where $\conj{\Phi}$ is the complex conjugate representation. Consider a $\bbZ$-lattice $\calM \subset D$. Scaling $\psi$ by a rational number, which does not change its unitary group, we can assume that
\begin{equation*}%\label{eq:PELlattice}
\tr_{D / \bbQ}(\psi(\calM, \calM)) = \bbZ.
\end{equation*}
We call $\Om = (D, \Phi, \sig; \psi, \calM)$ a \emph{PEL-type}.

A \emph{PEL-structure} $(\calA, \calC, \theta)$ of type $\Om$ consists of an abelian variety $\calA$ with polarization $\calC$ and $\theta : D \to \End_\bbQ(\calA)$ satisfying the following properties:
\begin{itemize}\itemsep-0.25em

\item There is a complex torus $\bbC^n / \Lam$, an isomorphism $\eta : D_\bbR \to \bbC^n$, and a homomorphism $\iota : \bbC^n \to \calA$ so that $\iota$ induces a biregular isomorphism from $\bbC^n / \Lam$ to $\calA$ such that
\begin{equation*}%\label{eq:PELtheta}
\iota(\Phi(a) w) = \thet(a) \iota(w)
\end{equation*}
for all $a \in D$ and $w \in \bbC^n$.

\item For all $a \in D$ and $x \in D_\bbR$, one has $\eta(\calM) = \Lam$ and $\eta(a x) = \Theta(a) \eta(x)$.

\item The polarization $\calC$ contains a divisor that determines a Riemann form $E$ on $\bbC^n / \Lam$ such that
\begin{equation*}%\label{eq:PELpolar}
E(\eta(x), \eta(y)) = \mathrm{tr}_{D / \bbQ}(\psi(x,y))
\end{equation*}
for all $x, y \in D_\bbR$.

\end{itemize}
Any $\tau \in \Aut(\bbC)$ acts on PEL-structures by
\begin{equation*}%\label{eq:PELaut}
(\calA, \calC, \theta)^\tau = (\calA^\tau, \calC^\tau, \theta^\tau).
\end{equation*}
The PEL-type of $(\calA, \calC, \theta)^\tau$ is then $\Om^\tau = (D^\tau, \Phi^\tau, \sig^\tau; \psi^\tau, \calM^\tau)$.

Now, suppose that $\tau$ is the extension to $\Aut(\bbC)$ of an element of $\Aut(F)$ and that $\calM$ is in fact an order in $D$. We see that one can take $h^\prime = \tau(h)$. For any nonarchimedean place $\nu$ of $F$ and $\calO_\nu$ the integers of $F_\nu$,
\begin{equation}\label{eq:LocalOrder}
\calM_\nu = \calM \otimes_{\calO_F} \calO_\nu
\end{equation}
is a lattice in $D_\nu$, considered as a $F_\nu$ vector space, and the elements $P_\nu$ of $H_{F_\nu}$ stabilizing $\calM_\nu$ define an open compact subgroup. Considering $P_\nu$ as a subgroup of $H_{F_\nu} \le H(\bbA_F^\infty)$, its image in $H^\prime(\bbA_F^\infty)$ under the isomorphism from Theorem~\ref{thm:MilneSuh} must stabilize the lattice
\begin{equation}\label{eq:ConjLocLatt}
\calM^\tau_{\tau(\nu)} = \calM^\tau \otimes_{\calO_F} \calO_{\tau(\nu)},
\end{equation}
hence it maps to the corresponding open compact subgroup $P_{\tau(\nu)}$ of ${H^\prime}_{\hspace{-0.2em}F_{\tau(\nu)}}$ stabilizing $\calM^\tau_{\tau(\nu)}$. This proves:

%%%%%%%%%%%%%%%%%%%%
\begin{prop}\label{prop:MyForm}
With the notation established in this section, suppose that $\tau$ is the extension to $\Aut(\bbC)$ of an element in $\Aut(F)$. Then, up to similarity, one can take $h^\prime = \tau(h)$. Moreover, for any order $\calM$ of $D$, the isomorphism $H(\bbA_F^\infty) \to H^\prime(\bbA_F^\infty)$ in Theorem~\ref{thm:MilneSuh} sends the open compact subgroup of $H_{F_\nu}$ stabilizing the local lattice $\calM_\nu$ defined in Equation~\eqref{eq:LocalOrder} to the open compact subgroup of ${H^\prime}_{\hspace{-0.2em}F_{\tau(\nu)}}$ stabilizing the local lattice $\calM^\tau_{\tau(\nu)}$ defined in Equation~\eqref{eq:ConjLocLatt}.
\end{prop}
%%%%%%%%%%%%%%%%%%%%

%%%%%%%%%%%%%%%%%%%%
\section{Buildings and their vertices}\label{sec:Buildings}
%%%%%%%%%%%%%%%%%%%%

This section continues with the notation established in \S\ref{sec:Conj}. Our proof that certain fake projective plane groups have isomorphic profinite completions requires very specific information about the maximal compact subgroups in the arithmetic definition of each group. This section briefly introduces the building $X_\nu$ associated with the groups $H_{F_\nu}$ over the nonarchimedean completions $F_\nu$ of $F$ and the open compact subgroups of $H_{F_\nu}$ stabilizing vertices of $X_\nu$.

The algebraic groups $H_{F_\nu}$ we will need to understand are locally isomorphic to either $\SL_3(F_\nu)$ or the special unitary group $\SU(h, E_\mu)$ of a hermitian form $h$ with respect to a quadratic extension $E_\mu / F_\nu$. We handle each case individually.

%%%%%%%%%%%%%%%%%%%%
\subsection{The building for $\SL_3(F_\nu)$}\label{ssec:SL3building}
%%%%%%%%%%%%%%%%%%%%

This case is quite classical, and we refer to either \cite[\S 2.8]{Tits} or \cite[\S 2]{CSbuild} for details. The vertices of $X_\nu$ are homothety classes of lattices $\calL \subset F_\nu^3$. We denote the homothety class of $\calL$ by $[\calL]$. Here, a lattice means an $\calO_\nu$ submodule of rank three, where $\calO_\nu$ is the valuation ring of $F_\nu$. Choosing a uniformizer $\pi$ for $\calO_\nu$, distinct vertices $[\calL_1]$ and $[\calL_2]$ are adjacent if there are representatives $\calL_1$ and $\calL_2$ for the homothety classes so that $\pi \calL_1 \subset \calL_2 \subset \calL_1$. Chambers of $X_\nu$ consist of triples of mutually adjacent vertices.

Then $\PGL_3(F_\nu)$ acts transitively on the vertices of $X_\nu$. This implies that all vertices are \emph{hyperspecial}, and the stabilizer $P_\calL < \SL_3(F_\nu)$ of $[\calL]$ is a maximal open compact subgroup. We call $P_\calL$ a \emph{hyperspecial parahoric subgroup}, and note that all maximal parahoric subgroups are conjugate in $\GL_3(F_\nu)$. For consistency with the unitary case, we say that all vertices of the building are \emph{Type 1} vertices.

%%%%%%%%%%%%%%%%%%%%
\subsection{The building for $\SU(h, E_\mu)$}\label{ssec:SUbuilding}
%%%%%%%%%%%%%%%%%%%%

See \cite[\S 2]{CSbuild} for an excellent detailed construction of the building from our point of view. In this case, the building $X_\nu$ is a tree. As with $\SL_3(F_\nu)$, the vertices are related to homothety classes of lattices. Given a lattice $\calL$, let $\calL^*$ denote its dual with respect to the form $h$. It ends up that either $[\calL^*] = [\calL]$ or the two homothety classes define adjacent vertices in the building for $\SL_3(E_\mu)$. The vertices of $X_\nu$ fall into two classes:
%\begin{itemize}\itemsep-0.25em

\bigskip

\noindent
\emph{Type 1:} Classes $[\calL]$ for which $[\calL] = [\calL^*]$.

\bigskip

\noindent
\emph{Type 2:} Unordered pairs $([\calL], [\calL^*])$ of neighbors in the building for $\SL_3(E_\mu)$.

\bigskip

%\end{itemize}
A Type 1 vertex $[\calL^\prime]$ is adjacent to a Type 2 vertex $([\calL], [\calL^*])$ if and only if the triple forms a chamber in the building for $\SL_3(E_\mu)$. A vertex is never adjacent to another of the same type. Type 1 vertices are hyperspecial and the Type 2 vertices are special, but not hyperspecial. The action of $\SU(h, E_\mu)$ preserves the type of a vertex, and the action is transitive on vertices of a given type. The maximal open compact subgroups of $\SU(h, E_\mu)$ are vertex stabilizers, hence there are two conjugacy classes.

%%%%%%%%%%%%%%%%%%%%
\section{Isomorphic algebraic fundamental groups}\label{sec:ThePairs}
%%%%%%%%%%%%%%%%%%%%

We now prove the isomorphisms between algebraic fundamental groups that, combined with Theorem~\ref{thm:SameList}, complete the proof of Theorem~\ref{thm:Main}. The general progression for each pair $\{j,k\}$ is:
\begin{enum}\itemsep-0.25em

\item Describe the number field $F$ and $F$-algebraic group $H$ associated with the pair.

\item Relate $\Gam_j$ and $\Gam_k$ to certain congruence subgroups of $H(F)$.

\item Study the action of the extension to $\Aut(\bbC)$ of a generator $\tau$ for $\Gal(F / \bbQ)$ on the congruence quotients of the ball related to $\Gam_j$ and $\Gam_k$ using Theorem~\ref{thm:MilneSuh}.

\item When $\Gam_j$ and $\Gam_k$ do not themselves arise from congruence subgroups of $H(F)$, use additional covering space arguments to study the action of $\tau$ on the relevant fake projective planes.

\end{enum}
It follows in each case that $\wh{\Gam}_j \cong \wh{\Gam}_k$, as desired. When Theorem~\ref{thm:MilneSuh} does not apply directly to the lattices at hand we must use nearby auxiliary lattices, which introduces some technical annoyances. In an attempt to make each case relatively self-contained, there is some near-verbatim repetition in the presentation of each case. However each case also requires at least one ingredient not seen in the other three.

%%%%%%%%%%%%%%%%%%%%
\subsection{The pair $\{34, 35\}$}
%%%%%%%%%%%%%%%%%%%%

We follow \cite[\S 9.3]{PrasadYeung}. See \cite{CartwrightC2, CartwrightpadicC2} for a more concrete description of these groups as matrix groups. This is the most complicated pair of the four, so a first-time reader may want to read the later cases first.

%%%%%%%%%%%%%%%%%%%%
\subsubsection{The number fields and some prime ideals}\label{ssec:C2fields}
%%%%%%%%%%%%%%%%%%%%

Let $F = \bbQ(\al)$ with $\al^2 = 5$ and $E = F(\beta)$ with $\beta^2 = -3$. There is a unique prime ideal $\frakp_2$ of $F$ dividing $2$, and $\frakp_2$ splits as a product of two primes $\frakq_{2+}$ and $\frakq_{2-}$ of $E$. All three primes have residue field $\bbF_4$.

%%%%%%%%%%%%%%%%%%%%
\subsubsection{The algebraic group and open compact subgroups}\label{ssec:C2gp}
%%%%%%%%%%%%%%%%%%%%

The simply connected $F$-algebraic group $H$ is defined in \cite[\S 9.1]{PrasadYeung} using the unique cubic division algebra $D$ over $E$ with nontrivial invariants at the two primes of $E$ dividing $2$ and a hermitian element $h \in D$. We refer to \cite{PrasadYeung} for details of the construction, which will not concern us. We now define open compact subgroups of $H_{F_\nu}$ for each nonarchimedean place $\nu$.

If $|\nu|$ denotes the residue characteristic of $\nu$ and $|\nu| \ge 3$, we take $P_\nu$ to be a hyperspecial parahoric subgroup, that is, the stabilizer of a Type 1 vertex of $X_\nu$. It remains to define the subgroup $P_2 < H_{F_2}$, where $F_2$ is the completion associated with $\frakp_2$. For each prime $\frakq_{2 \pm}$ of $E$ dividing $\frakp_2$, the completion $D_{2 \pm}$ of $D$ is a division algebra. The units of norm one in $D_{2\pm}$ are isomorphic as abstract groups and are mutually isomorphic to $H_{F_2}$. Specifically,
\begin{equation*}%\label{eq:C2gp2}
H_{F_2} \cong \left\{(x, y) \in D_{2+}^1 \times D_{2-}^1\ :\ \sig(x, y) (h, h) (x, y) = (h, h)\right\},
\end{equation*}
where $\sig$ extends to the factor-swapping isomorphism, so $\sig$-hermitian elements of $D$ map in diagonally. It follows that $y = h x h^{-1}$ defines an isomorphism between $H_{F_2}$ and $D_{2+}^1$. Under the above isomorphism of $H_{F_2}$ with $D_{2+}^1$, let $\calO_{2+}$ be the unique maximal order of $D_{2+}$. Then $P_2 \cong \calO_{2+}^1$.

%%%%%%%%%%%%%%%%%%%%
\subsubsection{The arithmetic groups}\label{ssec:C2lattices}
%%%%%%%%%%%%%%%%%%%%

Let $K = \prod P_\nu < H(\bbA_F^\infty)$ with $P_\nu$ defined in \S\ref{ssec:C2gp} and set $\Gam = H(F) \cap K$. Let $\Lam < \SU(2,1)$ be the normalizer of $\Gam$ and $\conj{\Gam}, \conj{\Lam}$ denote the images of these groups in $\PU(2,1)$. The group $\Gam_{32}$ is the index $3$ normal congruence subgroup of $\conj{\Gam}$ described in \cite[Rem.~9.4]{PrasadYeung}. In particular, $\Gam_{32}$ is the congruence subgroup of $H^{ad}(F)^+$ associated with the open compact subgroup
\begin{equation}\label{eq:32oc}
K^\prime = P_2^\prime \times \prod_{|\nu| \ge 3} P_\nu < H(\bbA_F^\infty),
\end{equation}
where $P_2^\prime$ is the kernel of the unique homomorphism from $P_2$ to $\bbZ / 3$.

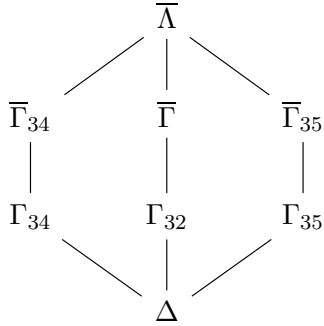
\begin{figure}[h]
\centering
\begin{tikzcd}
& \conj{\Lam} \arrow[-, d] \arrow[-, dr] \arrow[-, dl] & \\
\conj{\Gam}_{34} \arrow[-, d] & \conj{\Gam} \arrow[-, d] & \conj{\Gam}_{35} \arrow[-, d] \\
\Gam_{34} \arrow[-, dr] & \Gam_{32} \arrow[-, d] & \Gam_{35} \arrow[-, dl] \\
& \Del & 
\end{tikzcd}
\caption{The diagram of fake projective plane groups for $j = 32, 34, 35$}\label{fig:32-35}
\end{figure}

The groups $\Gam_j$ for $j = 32, 34, 35$ then fit into the diagram given in Figure~\ref{fig:32-35}. All neighboring inclusions in the diagram are index $3$. We now give a sequence of facts that combine information from \cite[\S 9.3]{PrasadYeung}, Cartwright and Steger's presentation for $\conj{\Lam}$ \cite{CartwrightgensC2}, and Magma to give a more detailed description of the groups in Figure~\ref{fig:32-35}.

%%%%%%%%%%%%%%%%%%%%
\begin{fact}\label{C2facta}
The group $\conj{\Lam}$ contains a unique subgroup of index $3$ admitting a homomorphism onto $\bbZ / 21$. This must be $\conj{\Gam}$.
\end{fact}
%%%%%%%%%%%%%%%%%%%%

%%%%%%%%%%%%%%%%%%%%
\begin{fact}\label{C2factb}
The group $\conj{\Gam}$ contains a unique normal subgroup of index $3$. This must be $\Gam_{32}$.
\end{fact}
%%%%%%%%%%%%%%%%%%%%

%%%%%%%%%%%%%%%%%%%%
\begin{fact}\label{C2factc}
The group $\Del = \Gam_{32} \cap \Gam_{34}$ equals $\Gam_{32} \cap \Gam_{35}$. Then $\Del$ is a normal index $3$ subgroup of both $\Gam_{34}$ and $\Gam_{35}$, but is not normal in $\Gam_{23}$. Among all conjugacy classes of index $3$ subgroups of $\Gam_{32}$, $\Del$ is in the unique conjugacy class that admits four distinct homomorphisms onto $S_3$.
\end{fact}
%%%%%%%%%%%%%%%%%%%%

%%%%%%%%%%%%%%%%%%%%
\begin{fact}\label{C2factd}
The lift of any element of $\conj{\Gam}$ to $D^*$ has reduced norm contained in $(E^*)^3$.
\end{fact}
%%%%%%%%%%%%%%%%%%%%

%%%%%%%%%%%%%%%%%%%%
\begin{fact}\label{C2facte}
In terms of Cartwright and Steger's generators $\mathrm{AA}$, $\mathrm{BB}$, $\mathrm{ZZ}$ for $\conj{\Gam}$ in \cite{CartwrightgensC2}, coset representatives for $\Del$ in $\Gam_{34}$ are given by $(\mathrm{BB}\, \mathrm{ZZ}^{-1})^{\pm 2}$.
\end{fact}
%%%%%%%%%%%%%%%%%%%%

%%%%%%%%%%%%%%%%%%%%
\begin{fact}\label{C2factf}
In terms of Cartwright and Steger's generators $\mathrm{AA}$, $\mathrm{BB}$, $\mathrm{ZZ}$ for $\conj{\Gam}$ in \cite{CartwrightgensC2}, coset representatives for $\Del$ in $\Gam_{35}$ are given by $(\mathrm{BB}\, \mathrm{ZZ})^{\pm 2}$.
\end{fact}
%%%%%%%%%%%%%%%%%%%%

Let $\xi = \frac{1}{2}(-1 + \beta)$ be a primitive $3^{rd}$ root of unity and set
\begin{equation*}
b = \det(\mathrm{BB}) = \frac{1 - \al \beta}{4}.
\end{equation*}
By Fact~\ref{C2factd}, any lift to $D^*$ of an element of $\Del$ has reduced norm contained in $(E^*)^3$. Then $\det(\mathrm{Z}) = \xi^2$, and Facts~\ref{C2facte}-\ref{C2factf} imply that $\Gam_{34}$ is distinguished from $\Gam_{35}$ by the property that every preimage $\wh{\gam}$ of $\gam \in \Gam_j$ in $D^*$ has reduced norm contained in the subgroup
\begin{equation*}%\label{eq:C18detcond}
N_j = \left\langle b^{2 m_j} \xi (E^*)^3 \right\rangle \subset E^*
\end{equation*}
where $m_{34} = 1$ and $m_{35} = 2$. The identification of $\Gam_j / \Del$ with $\bbZ / 3$ is given by the isomorphism $N_j / (E^*)^3 \cong \bbZ / 3$.

%%%%%%%%%%%%%%%%%%%%
\subsubsection{The isomorphism of profinite completions for $j = 34, 35$}
%%%%%%%%%%%%%%%%%%%%

%%%%%%%%%%%%%%%%%%%%
\begin{prop}\label{prop:34n35}
Let $X_{34}, X_{35}$ be fake projective planes with fundamental groups $\Gam_{34}, \Gam_{35}$ and choose a representative $\tau \in \Aut(\bbC)$ for the generator of $\Gal(F / \bbQ)$. Then, up to complex conjugation, $X_{34}^\tau$ is isomorphic to $X_{35}$. Therefore $\wh{\Gam}_{34} \cong \wh{\Gam}_{35}$.
\end{prop}
%%%%%%%%%%%%%%%%%%%%

%%%%%%%%%%%%%%%%%%%%
\begin{pf}
Choosing the appropriate complex structures, the diagram of subgroups in Figure~\ref{fig:32-35} implies that we have the diagram of $\bbZ / 3$ \'etale coverings in Figure~\ref{fig:2maps}, where $Y$ is the ball quotient associated with $\Del$. Let $p_j : Y \to X_j$ be the covering and $\psi_j$ a generator for its Deck group. We will prove, as indicated in Figure~\ref{fig:2maps}, that $Y^\tau \cong Y$, which implies that $\tau$ acts on $\Aut(Y)$. We then show that, up to the correct choice of generators for the Deck groups, $\psi_{34} ^\tau = \psi_{35}$ in the sense of \S\ref{sec:Action}. This implies that $p_{34}^\tau = p_{35}$, hence $X_{34}^\tau \cong X_{35}$, which proves the proposition.
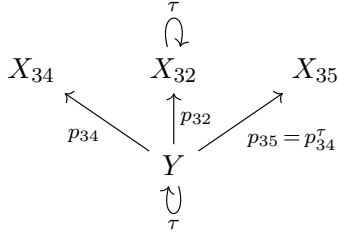
\begin{figure}[h]
\centering
\begin{tikzcd}
& X_{34} & X_{32} \arrow[loop above, "\tau"] & X_{35} & \\
& & Y \arrow[ul, "p_{34}"] \arrow[ur, "p_{35}\, =\, p_{34}^\tau" below right] \arrow[loop below, "\tau"] \arrow[u, "p_{32}" right] & &
\end{tikzcd}
\caption{The common $\bbZ / 3$ cover of $X_{32}$, $X_{34}$, and $X_{35}$}\label{fig:2maps}
\end{figure}

To prove that $Y^\tau \cong Y$, we first show that $X_{32}^\tau \cong X_{32}$, which is a consequence of Proposition~\ref{prop:MyForm}. To see this, $\Gam_{32}$ is the group associated with the open compact subgroup $K^\prime$ of $H(\bbA_F^\infty)$ defined in Equation~\eqref{eq:32oc}. For places $\nu$ with $|\nu| \ge 3$, $P_\nu$ is the stabilizer of a Type 1 vertex of the building $X_\nu$ (i.e., $P_\nu$ is hyperspecial), hence the image of $K$ under the self-isomorphism of $H(\bbA_F^\infty)$ induced by $\tau$ is hyperspecial at all such $\nu$. For $\frakp_2$, the image of $P_2$ is the unique maximal order of $D_{2+} \cong D_{2-}$, so it is again $P_2$ under the isomorphism induced by $\tau$. Since $P_2^\prime$ is the kernel of the unique homomorphism from $P_2$ to $\bbZ / 3$, the image of $P_2^\prime$ is also isomorphic to $P_2^\prime$. This implies that, up to conjugacy in $H(\bbA_F^\infty)$, the self-isomorphism takes $K^\prime$ to $K^\prime$. In other words, $X_{32}^\tau$ has fundamental group $\Gam_{32}$ and so $X_{32}^\tau \cong X_{32}$.

Then, applying Fact~\ref{C2factc} to $Y^\tau$, it is the unique connected \'etale cover of $X_{32}^\tau = X_{32}$ whose fundamental group admits four distinct homomorphisms onto $S_3$. By Theorem~\ref{thm:QuoNum} we see that, up to conjugacy in $\Gam_{32}$, $\pi_1(Y^\tau)$ equals $\pi_1(Y)$. Therefore, $Y^\tau \cong Y$, which proves the first assertion.

We now study $\psi_{34}^\tau$ by lifting to a self-isometry of the universal cover. A choice of lift of $\psi_{34}$ is given by choosing a coset representative $\gam_{34}$ for $\conj{\Gam}$, where the choice of nontrivial coset representative is equivalent to a choice of generator for the Deck group. Following the notation in \S\ref{ssec:C2lattices}, we choose $\gam_{34}$ so that any lift of $\gam_{34}$ to $D^*$ has determinant contained in $b^2 \xi (E^*)^3$.

As in \cite[\S 5.2]{Miyake}, we identify complex hyperbolic space $\bbB^2$ with the space $\{\calQ_z\}_{z \in \bbB^2}$ of PEL-structures of type $\Om$, with $\Om$ defined from the algebra $D$ and hermitian form $h$ as in \S\ref{sec:Conj}. The lattice $\calM$ associated with the PEL-type $\Om$ is stabilized by $\Gam \le D^1$. The lift of $\psi_{34}^\tau$ to the universal cover of $X_{34}^\tau$ is then given by:
\begin{align*}
\calQ_z^\tau &\mapsto \left(\tau \circ \gam_{34} \circ \tau^{-1}\right) \calQ_z^\tau \\ %\label{eq:43taulift} \\
&\, = \gam_{34}^\tau \calQ_z^\tau, %\nonumber
\end{align*}
where $\gam_{34}^\tau$ means in the sense of the $\Aut(\bbC)$ action as in Equation~\eqref{eq:Marking}, not the $\Gal(F / \bbQ)$ action on $D$.

By Corollary~\ref{cor:FakeQuotient}, $\gam_{34}^\tau$ is the lift to $\bbB^2$ of an element $\Aut(Y)$ of order $3$ with no fixed points and quotient a fake projective plane. We can then choose a lift $\wh{\gam}_{34}$ of $\gam_{34}$ to $D^*$, and the action of $\gam_{34}^\tau$ on PEL-types takes ${\Om^\tau = (D^\tau, \Phi^\tau, \sig^\tau; \psi^\tau, \calM^\tau)}$ to the equivalent PEL-type
\begin{equation*}\Om^\prime = (D^\tau , \Phi^\tau, \sig^\tau; \psi^\tau, (\calM \wh{\gam}_{34})^\tau) %\label{eq:C2congPEL}
\end{equation*}
(see the Corollary in \S 5.4 of \cite{Miyake}), where now we can consider $(\calM \wh{\gam}_{34})^\tau$ as the $\tau$-conjugate for the Galois action on $\calM \subset D$ by identifying $D_\bbR$ with $\bbC^9$. Taking the determinant, this implies that a lift $\wh{\gam}_{34}^\tau$ of $\gam_{34}^\tau$ to $D^*$ has reduced norm contained in $\tau(\det(\wh{\gam}_{34})) (E^*)^3 = \tau(b^2 \xi) (E^*)^3$.

Then $\tau(b^2) = b^{-2}$ and $\tau(\xi) = \xi$, so any representative for $\gam_{34}^\tau$ in $D^*$ has reduced norm contained in $b \xi (E^*)^3$. This is a nontrivial coset representative for $\Del$ in $\Gam_{35}$ by Fact~\ref{C2factf}, hence (up to the correct choice of generator) $\psi_{34}^\tau = \psi_{35}$. This completes the second and final step required to prove the proposition.
\end{pf}
%%%%%%%%%%%%%%%%%%%%

%%%%%%%%%%%%%%%%%%%%%
%\begin{prop}\label{prop:34n35}
%Let $X_{34}, X_{35}$ be fake projective planes with fundamental groups $\Gam_{34}, \Gam_{35}$ and choose a representative $\tau \in \Aut(\bbC)$ for the generator of $\Gal(F / \bbQ)$. Then, up to complex conjugation, $X_{34}^\tau$ is isomorphic to $X_{35}$. Therefore $\wh{\Gam}_{34} \cong \wh{\Gam}_{35}$.
%\end{prop}
%%%%%%%%%%%%%%%%%%%%%
%
%%%%%%%%%%%%%%%%%%%%%
%\begin{pf}
%Equations for these fake projective planes were computed by Borisov and Fatighenti. See \cite{BF} and the auxiliary files referenced there. Their equations for the pair are equivalent under the nontrivial Galois involution of $F / \bbQ$. Taking an extension of this element to $\Aut(\bbC)$ we see that, up to the correct choice of complex structure on each, $X_{34}$ and $X_{35}$ are $\Aut(\bbC)$ equivalent. This proves the proposition.
%\end{pf}
%%%%%%%%%%%%%%%%%%%%%

%%%%%%%%%%%%%%%%%%%%
\begin{rem}\label{rem:BF}
Proposition~\ref{prop:34n35} should also follow from work on Borisov and Fatighenti. See \cite{BF} and the auxiliary files referenced there. Their work gives equations for the relevant pair of fake projective planes that are invariant under $\Gal(F/\bbQ)$. It seems difficult to rigorously verify that the $\Gal(F / \bbQ)$-orbit of equations derived from \cite{BF} are indeed for two distinct fake projective planes, so we give the above proof instead.
\end{rem}
%%%%%%%%%%%%%%%%%%%%

%%%%%%%%%%%%%%%%%%%%
\subsection{The pair $\{43, 44\}$}
%%%%%%%%%%%%%%%%%%%%

This case is arguably the easiest of the four. We follow \cite[\S A.6]{PYFix}. See \cite{CartwrightC18, CartwrightpadicC18} for a more concrete description of these groups as matrix groups.

%%%%%%%%%%%%%%%%%%%%
\subsubsection{The number fields and some prime ideals}\label{ssec:C18fields}
%%%%%%%%%%%%%%%%%%%%

Let $F = \bbQ(\al)$ with $\al^2 = 6$ and $E = F(\beta)$ with $\beta^2 = -3$. There is a unique prime ideal $\frakp_2$ of $F$ dividing $2$, and $\frakp_2$ is divisible by a unique prime $\frakq_2$ of $E$. The residue fields are $\bbF_2$ and $\bbF_4$, respectively. There is also a unique prime $\frakp_3$ of $F$ dividing $3$, which is ramified over $\bbQ$ with residue field $\bbF_3$. Then $\frakp_3$ splits into two primes $\frakq_{3 \pm}$ of $E$, which also have residue field $\bbF_3$.

%%%%%%%%%%%%%%%%%%%%
\subsubsection{The algebraic group and open compact subgroups}\label{ssec:C18gp}
%%%%%%%%%%%%%%%%%%%%

The simply connected $F$-algebraic group $H$ is defined in \cite[\S 9.1]{PrasadYeung} using the unique cubic division algebra $D$ over $E$ with nontrivial invariants at the two primes of $E$ dividing $3$ and a hermitian element $h \in D$. We refer to \cite{PrasadYeung} for details of the construction, which will not concern us. We now define open compact subgroups of $H_{F_\nu}$ for each nonarchimedean place $\nu$.

Let $F_2$ be the completion of $F$ associated with $\frakp_2$. Then $H_{F_2}$ is the special unitary group of $h$ with respect to the extension $E_2 / F_2$, where $E_2$ is the associated completion of $E$. We take $P_2$ to be the nonhyperspecial maximal parahoric subgroup stabilizing a Type 2 vertex of the building $X_2$ for $H_{F_2}$. If $|\nu|$ denotes the residue characteristic of $\nu$ and $|\nu| \ge 5$, we take $P_\nu$ to be a hyperspecial parahoric subgroup, that is, the stabilizer of a Type 1 vertex of $X_\nu$.

It remains to define the subgroup $P_3 < H_{F_3}$, where $F_3$ is the completion associated with $\frakp_3$. For each prime $\frakq_{3 \pm}$ of $E$ dividing $\frakp_3$, the completion $D_{3 \pm}$ of $D$ is a division algebra. The units of norm one in $D_{3\pm}$ are isomorphic as abstract groups and are mutually isomorphic to $H_{F_3}$. Then
\begin{equation*}%\label{eq:C18gp3}
H_{F_3} \cong \left\{(x, y) \in D_{3+}^1 \times D_{3-}^1\ :\ \sig(x, y) (h, h) (x, y) = (h, h)\right\},
\end{equation*}
where $\sig$ extends to the factor-swapping isomorphism, so $\sig$-hermitian elements of $D$ map in diagonally. It follows that $y = h x h^{-1}$ defines an isomorphism between $H_{F_3}$ and $D_{3+}^1$. Under the above isomorphism of $H_{F_3}$ with $D_{3+}^1$, let $\calO_{3+}$ be the unique maximal order of $D_{3+}$. Then $P_3 \cong \calO_{3+}^1$.

%%%%%%%%%%%%%%%%%%%%
\subsubsection{The arithmetic groups}\label{ssec:C18lattices}
%%%%%%%%%%%%%%%%%%%%

\begin{figure}[h]
\centering
\begin{tikzcd}
& \conj{\Lam} \arrow[-, d] \arrow[-, dr] \arrow[-, dl] & \\
\Gam_{43} \arrow[-, dr] & \Gam_{42} \arrow[-, d] & \Gam_{44} \arrow[-, dl] \\
& \conj{\Gam} &
\end{tikzcd}
\caption{The diagram of fake projective plane groups for $j = 42, 43, 44$}\label{fig:42-44}
\end{figure}
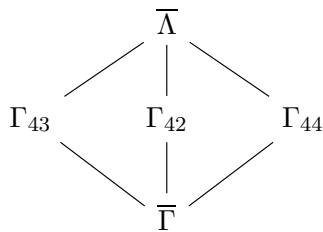

Let $K = \prod P_\nu \le H(\bbA_F^\infty)$ with $P_\nu$ defined in \S\ref{ssec:C18gp} and set $\Gam = H(F) \cap K$. Let $\Lam < \SU(2,1)$ be the normalizer of $\Gam$ and $\conj{\Gam}, \conj{\Lam}$ denote the images of these groups in $\PU(2,1)$. The groups $\Gam_j$ for $j = 42, 43, 44$ then fit into the diagram given in Figure~\ref{fig:42-44}, where all subgroups are normal. In fact, $\conj{\Lam} / \Gam_j$ and $\Gam_j / \conj{\Gam}$ are isomorphic to $\bbZ / 3$ for all $j$.

Let $\xi = \frac{1}{2}(-1 + \beta)$ be a primitive $3^{rd}$ root of unity and define
\begin{equation*}\label{eq:C18a}
a = (1 + \al \beta) (1 - \al \beta)^{-1}.
\end{equation*}
Then, as described in \cite[\S A.6]{PYFix}, $\Gam_j$ is characterized by the property that every preimage $\wh{\gam}$ of $\gam \in \Gam_j$ in $D^*$ has reduced norm contained in the subgroup 
\begin{equation*}%\label{eq:C18detcond}
N_j = \left\langle a^{m_j} \xi (E^*)^3 \right\rangle \subset E^*
\end{equation*}
for some $m_j \in \{0,1,2\}$. The identification of $\Gam_j / \conj{\Gam}$ with $\bbZ / 3$ is given by the isomorphism $N_j / (E^*)^3 \cong \bbZ / 3$. The notation in \cite{Cartwright} (see Table~\ref{tb:C18}) means that $m_{42} = 0$, $m_{43} = 1$, and $m_{44} = 2$.

%%%%%%%%%%%%%%%%%%%%
\subsubsection{The isomorphism of profinite completions for $j = 43, 44$}
%%%%%%%%%%%%%%%%%%%%

%%%%%%%%%%%%%%%%%%%%
\begin{prop}\label{prop:43n44}
Let $X_{43}, X_{44}$ be fake projective planes with fundamental groups $\Gam_{43}, \Gam_{44}$ and choose a representative $\tau \in \Aut(\bbC)$ for the generator of $\Gal(F / \bbQ)$. Then, up to complex conjugation, $X_{43}^\tau$ is isomorphic to $X_{44}$. Therefore $\wh{\Gam}_{43} \cong \wh{\Gam}_{44}$.
\end{prop}
%%%%%%%%%%%%%%%%%%%%

%%%%%%%%%%%%%%%%%%%%
\begin{pf}
Choosing the appropriate complex structures, the diagram of subgroups in Figure~\ref{fig:42-44} implies that we have the diagram of $\bbZ / 3$ \'etale coverings in Figure~\ref{fig:18maps}, where $Y$ is the ball quotient associated with $\conj{\Gam}$. Let $p_j : Y \to X_j$ be the covering and $\psi_j$ a generator for its Deck group. We will prove, as indicated in Figure~\ref{fig:18maps}, that $Y^\tau \cong Y$, which implies that $\tau$ acts on $\Aut(Y)$. We then show that, up to the correct choice of generators for the Deck groups, $\psi_{43} ^\tau = \psi_{44}$ in the sense of \S\ref{sec:Action}. This implies that $p_{43}^\tau = p_{44}$, hence $X_{43}^\tau \cong X_{44}$, which proves the proposition.
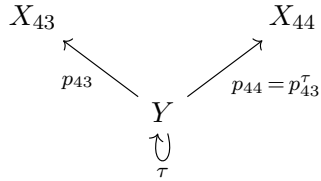
\begin{figure}[h]
\centering
\begin{tikzcd}
& X_{43} & & X_{44} & \\
& & Y \arrow[ul, "p_{43}"] \arrow[ur, "p_{44}\, =\, p_{43}^\tau" below right] \arrow[loop below, "\tau"] & &
\end{tikzcd}
\caption{The common $\bbZ / 3$ cover of $X_{43}$ and $X_{44}$}\label{fig:18maps}
\end{figure}

The fact that $Y^\tau \cong Y$ is a consequence of Proposition~\ref{prop:MyForm}. To see this, $\conj{\Gam} = \pi_1(Y)$ is the group associated with the open compact subgroup $K$ of $H(\bbA_F^\infty)$ defined in \S\ref{ssec:C18lattices}. For places $\nu$ with $|\nu| \ge 5$, $P_\nu$ is the stabilizer of a Type 1 vertex of the building $X_\nu$ (i.e., $P_\nu$ is hyperspecial), hence the image of $K$ under the self-isomorphism of $H(\bbA_F^\infty)$ induced by $\tau$ is hyperspecial at every $\nu$. For $\frakp_3$, the image open compact subgroup is defined by the unique maximal order of $D_{3+} \cong D_{3-}$, so it is again $P_3$ under the isomorphism induced by $\tau$. Lastly, $P_2$ is the stabilizer of a Type 2 vertex of the associated building, hence so is the image group. This implies that, up to conjugacy in $H(\bbA_F^\infty)$, the self-isomorphism takes $K$ to $K$. In other words, $Y^\tau$ has fundamental group $\conj{\Gam}$ and so $Y^\tau \cong Y$. This proves the first assertion.

We now study $\psi_{43}^\tau$ by lifting to a self-isometry of the universal cover. A choice of lift of $\psi_{43}$ is given by choosing a coset representative $\gam_{43}$ for $\conj{\Gam}$ in $\Gam_{43}$, where the choice of nontrivial coset representative is equivalent to a choice of generator for the Deck group. Following the notation in \S\ref{ssec:C18lattices}, we choose $\gam_{43}$ so that any lift of $\gam_{43}$ to $D^*$ has determinant contained in $a \xi (E^*)^3$.

As in \cite[\S 5.2]{Miyake}, we identify complex hyperbolic space $\bbB^2$ with the space $\{\calQ_z\}_{z \in \bbB^2}$ of PEL-structures of type $\Om$, with $\Om$ defined from the algebra $D$ and hermitian form $h$ as in \S\ref{sec:Conj}. The lattice $\calM$ associated with the PEL-type $\Om$ is stabilized by $\Gam \le D^1$. The lift of $\psi_{43}^\tau$ to the universal cover of $Y_{43}^\tau$ is then given by:
\begin{align*}
\calQ_z^\tau &\mapsto \left(\tau \circ \gam_{43} \circ \tau^{-1}\right) \calQ_z^\tau \\ %\label{eq:43taulift} \\
&\, = \gam_{43}^\tau \calQ_z^\tau, %\nonumber
\end{align*}
where $\gam_{43}^\tau$ means in the sense of the $\Aut(\bbC)$ action as in Equation~\eqref{eq:Marking}, not the $\Gal(F / \bbQ)$ action on $D$.

By Corollary~\ref{cor:FakeQuotient}, $\gam_{43}^\tau$ is the lift to $\bbB^2$ of an element $\Aut(Y)$ of order $3$ with no fixed points and quotient a fake projective plane. We can then choose a lift $\wh{\gam}_{43}$ of $\gam_{43}$ to $D^*$, and the action of $\gam_{34}^\tau$ on PEL-types takes ${\Om^\tau = (D^\tau, \Phi^\tau, \sig^\tau; \psi^\tau, \calM^\tau)}$ to the equivalent PEL-type
\begin{equation*}\Om^\prime = (D^\tau , \Phi^\tau, \sig^\tau; \psi^\tau, (\calM \wh{\gam}_{43})^\tau) %\label{eq:C2congPEL}
\end{equation*}
(see the Corollary in \S 5.4 of \cite{Miyake}), where now we can consider $(\calM \wh{\gam}_{43})^\tau$ as the $\tau$-conjugate for the Galois action on $\calM \subset D$ by identifying $D_\bbR$ with $\bbC^9$. Taking the determinant, this implies that a lift $\wh{\gam}_{43}^\tau$ of $\gam_{43}^\tau$ to $D^*$ has reduced norm contained in $\tau(\det(\wh{\gam}_{43})) (E^*)^3 = \tau(a \xi) (E^*)^3$.

Then $\tau(a) = a^{-1}$ and $\tau(\xi) = \xi$, so any representative for $\gam_{43}^\tau$ in $D^*$ has reduced norm contained in $a^2 \xi (E^*)^3$. This is a nontrivial coset representative for $\conj{\Gam}$ in $\Gam_{44}$, hence (up to the correct choice of generator) $\psi_{43}^\tau = \psi_{44}$. This completes the second and final step required to prove the proposition.
\end{pf}
%%%%%%%%%%%%%%%%%%%%

%%%%%%%%%%%%%%%%%%%%
\subsection{The pairs $\{47, 49\}$ and $\{48, 50\}$}
%%%%%%%%%%%%%%%%%%%%

These pairs are commensurable and best handled together. The argument for these cases benefits from slightly more refined information about the commensurability class than was needed for previous cases. We follow \cite{CartwrightC20}, but change some notation.

%%%%%%%%%%%%%%%%%%%%
\subsubsection{The number fields}\label{ssec:C20fields}
%%%%%%%%%%%%%%%%%%%%

Let $F = \bbQ(\al)$ with $\al^2 = 7$ and $E = F(\beta)$ with $\beta^2 = -1$. If $\zeta$ is a primitive $7^{th}$ root of unity and $L = E(\zeta)$, then $L / E$ is a cubic Galois extension with Galois group generated by the element $\phi$ characterized by
\begin{equation*}%\label{eq:C20GalGp}
\phi(\zeta) = \zeta^2.
\end{equation*}
The maximal totally real subfield of $L$ is $L^+ = F(\om)$, where
\begin{equation*}%\label{eq:C20om}
\om = \zeta + \zeta^{-1},
\end{equation*}
which has minimal polynomial $t^3 + t^2 - 2 t - 1$. Then $\phi(\om) = \om^2 - 2$ and the restriction of $\phi$ to $L^+$ generates the Galois group of $L^+ / F$.

%%%%%%%%%%%%%%%%%%%%
\subsubsection{Some prime ideals}\label{ssec:C20primes}
%%%%%%%%%%%%%%%%%%%%

There is a unique prime ideal $\frakp_2$ of $F$ dividing $2$, and $\frakp_2$ splits as a product of two primes $\frakq_{2+}$ and $\frakq_{2-}$ of $E$. All three primes have residue field $\bbF_2$, and all three ideals area principal; we can take $\pi_2 = 3 + \al$ as a generator for $\frakp_2$, $\sig_{2+} = 1 + \frac{\al - \beta}{2}$ as a generator for $\frakq_{2+}$, and $\sig_{2-} = 1 - \frac{\al + \beta}{2}$ as a generator for $\frakq_{2-}$. Both of $\frakq_{2\pm}$ remain prime in $L$, hence $\sig_{2\pm}$ are generators of prime ideals of $L$ with residue field $\bbF_8$.

There are two primes $\frakp_{3\pm}$ of $F$ dividing $3$, and each defines a unique prime ideal $\frakq_{3\pm}$ of $E$. The residue fields are $\bbF_3$ and $\bbF_9$, respectively. One can take $\pi_{3\pm} = 2 \pm \al$ as a generator for either $\frakp_{3\pm}$ or $\frakq_{3\pm}$. This ordering is chosen so that $\al$ reduces modulo $\frakp_{3\pm}$ to $\pm 1$. The completion of $F$ with respect to either place $\nu_{3 \pm}$ associated with $\frakp_{3\pm}$ is isomorphic to $\bbQ_3$ and the corresponding completion of $E$ is isomorphic to the unique unramified quadratic extension $E_3$ of $\bbQ_3$. There are two primes of $L$ dividing $3$ with residue field $\bbF_{729}$, and the completion of $L$ with respect to either prime is the unique unramified cubic extension $L_3$ of $E_3$.

There is a unique prime $\frakp_7$ of $F$ dividing $7$, which has generator $\pi_7 = \al$ and residue field $\bbF_7$. Then $\pi_7$ continues to generate the unique prime ideal $\frakq_7$ of $E$ dividing $\frakp_7$, which has residue field $\bbF_{49}$. The unique prime ideal $\wh{\frakq}_7$ of $L$ dividing $\frakq_7$ has generator $\wh{\pi}_7 = \zeta - 1$ and residue field $\bbF_{49}$ (i.e., the prime $\frakq_7$ is totally ramified in $L$).

%%%%%%%%%%%%%%%%%%%%
\subsubsection{The algebra}\label{ssec:C20alg}
%%%%%%%%%%%%%%%%%%%%

Consider the cubic division algebra $D$ over $E$ generated by $L$ and an element $y$ subject to the relations:
\begin{align*}
y x y^{-1} &= \phi(x) \textrm{ for all } x \in L \\ %\label{eq:C20Rel1} \\
y^3 &= \frac{3 + \al \beta}{4} %\label{eq:C20Rel2}
\end{align*}
Set $\de = y^3$. Then
\begin{align}
\iota(x) &= \begin{pmatrix} x & 0 & 0 \\ 0 & \phi(x) & 0 \\ 0 & 0 & \phi^2(x) \end{pmatrix} \\ \label{eq:C20Emb1} \\
\iota(y) &= \begin{pmatrix} 0 & 1 & 0 \\ 0 & 0 & 1 \\ \de & 0 & 0 \end{pmatrix} \label{eq:C20Emb2}
\end{align}
is an embedding of $D$ in $\M_3(L)$ such that the subgroup $D^1$ of elements of reduced norm one are precisely the elements that map to $\SL_3(L)$. The assignments
\begin{align*}
\sig(\al) &= \al & \sig(\om) &= \om \\
\sig(\zeta) &= \zeta^{-1} & \sig(y) &= y^{-1}
\end{align*}
extend to an anti-involution of $D$ that is the restriction of conjugate transposition on $\M_3(L)$ with respect to $\Gal(L / L^+)$ under either embedding.

Using the factorization
\begin{equation}\label{eq:C20de}
\de = (3 \al - 8) \sig_{2+}^2 \sig_{2-}^{-2},
\end{equation}
with $3 \al - 8$ a unit of $F$, one checks that $D$ remains a division algebra over the localization of $E$ at each of $\frakq_{2\pm}$ and that these are the unique places of $E$ where $D$ remains a division algebra; see \cite{CartwrightC20}. Let $F_2$ be the completion of $F$ at the place associated with $\frakp_2$, and let $D_{2 \pm}$ be the completion of $D$ with respect to the completion of $E$ associated with $\frakq_{2\pm}$. Then $D_{2-}$ is the opposite algebra of $D_{2+}$. For every other place $\mu$ of $E$ with completion $E_\mu$, $D \otimes_E E_\mu$ is isomorphic to $\M_3(E_\mu)$.

%%%%%%%%%%%%%%%%%%%%
\begin{rem}\label{rem:C20DInt}
If one wants $y^3$ to be an algebraic integer, $\de$ could be replaced with $(3 \al - 8) \sig_{2+}^2 \sig_{2-} = \sig_{2-}^3 \de$.
\end{rem}
%%%%%%%%%%%%%%%%%%%%

%%%%%%%%%%%%%%%%%%%%
\subsubsection{The hermitian element and algebraic group}\label{ssec:C20gp}
%%%%%%%%%%%%%%%%%%%%

The element
\begin{align*}
h &= \frac{3 - \al - 3 \beta + \al \beta}{2} + \frac{3 - \al - 3 \beta + \al \beta}{2} \zeta + \left(1 - \frac{\al + \beta}{2}\right)\zeta^2 \\ %\label{eq:C20h} \\
&= \left(2-\frac{4 \al}{7}\right) + \left(1-\frac{3 \al}{7}\right)\om - \frac{\al}{7} \om^2 %\nonumber
\end{align*}
of $L^+ \subset D$ maps to a hermitian element of $\M_3(L)$ with respect to the Galois involution of $L / L^+$. Let $\nu_j : F \to \bbR$ denote the embedding generated by
\begin{equation*}%\label{C20:Emb3}
\nu_j(\al) = (-1)^j \sqrt{7}
\end{equation*}
for $j \in \{0,1\}$. Then $\nu_j$ extends to a complex conjugate pair of embeddings of $E$ into $\bbC$. Fix one such extension (e.g., $\nu_j(\beta) = i$). Defining $\nu_j(\om) = 2 \cos(2 \pi / 7)$ further extends $\nu_j$ to a complex embedding of $L$ (the other extensions are obtained by precomposition with $\phi$). Applying $\nu_j$ to matrices then gives us a pair of embeddings, still denoted by $\nu_j$, of $D$ into $\M_3(\bbC)$. One checks that $\nu_0(h)$ is a hermitian form of signature $(1,2)$ and $\nu_1(h)$ has signature $(3,0)$.

Our simply connected semisimple $F$-algebraic group $H$ is then the special unitary group of $h$. That is,
\begin{equation*}%\label{eq:C20Group}
H(F) = \left\{z \in D^1\ :\ \sig(z) h z = h\right\}
\end{equation*}
and $H_*(\bbR) \cong \SU(1,2) \times \SU(3)$. We now define the open compact subgroups $K_j$ of $H(\bbA_F^\infty)$ for $j \in \{47, 49\}$.

%%%%%%%%%%%%%%%%%%%%
\subsubsection{The open compact subgroup for $\frakp_2$}\label{ssec:C20oc2}
%%%%%%%%%%%%%%%%%%%%

The units of norm one in $D_{2\pm}$ are isomorphic as abstract groups and are mutually isomorphic to $H_{F_2}$. Then
\begin{equation*}%\label{eq:C20gp2}
H_{F_2} \cong \left\{(x, y) \in D_{2+}^1 \times D_{2-}^1\ :\ \sig(x, y) (h, h) (x, y) = (h, h)\right\},
\end{equation*}
where $\sig$ extends to the factor-swapping isomorphism, so $\sig$-hermitian elements of $D$ map in diagonally. It follows that $y = h x h^{-1}$ defines an isomorphism between $H_{F_2}$ and $D_{2+}^1$. Under the above isomorphism of $H_{F_2}$ with $D_{2+}^1$, let $\calO_{2+}$ be the unique maximal order of $D_{2+}$. Then $P_2 \cong \calO_{2+}^1$.

%%%%%%%%%%%%%%%%%%%%
\subsubsection{The open compact subgroups for $|\nu| \ge 5$}\label{ssec:C20ocother}
%%%%%%%%%%%%%%%%%%%%

Recall that $|\nu|$ denotes the residue characteristic of a place $\nu$ of $F$, and assume $|\nu| \ge 5$. For every place $\mu$ of $E$ dividing $\nu$ with completion $E_\mu$, $D \otimes_E E_\mu$ is isomorphic to $\M_3(E_\mu)$. For all places $\nu$ of $F$ other than $\frakp_2$, this implies
\begin{equation*}%\label{eq:C20Groups}
H_{F_\nu} \cong \begin{cases} \SU(h, E_\mu) & \nu \textrm{ inert in } E \\ \SL_3(F_\nu) & \nu \textrm{ splits in } E \end{cases}
\end{equation*}
where $\SU(h, E_\mu)$ denotes the special unitary group of $h$ with respect to the quadratic extension $E_\mu / F_\nu$. If $|\nu| \neq 7$, take $P_\nu$ to be the open compact subgroup consisting of all elements of $H_{F_\nu}$ that are $\nu$-adically integral with respect to the above basis. Then $P_\nu$ is hyperspecial for all such $\nu$.

Now, consider the unique place $\nu_7$ dividing $7$, which is associated with the prime $\frakp_7$ defined in \S\ref{ssec:C20primes}. There is a unique place $\mu_7$ of $E$ dividing $\nu_7$ and the associated local field extension $E_7 / F_7$ is an unramified quadratic extension. Note that the factorization in Equation~\eqref{eq:C20de} implies that $\de$ is a unit of $E_7$ (i.e., has valuation zero). Similarly, there is a unique place $\wh{\mu}_7$ of $L$ dividing $\mu_7$, and the completion $L_7$ is a totally ramified cubic extension of $E_7$. Even more, we will need the unique prime $\wh{\nu}_7$ of $L^+$ divisible by $\wh{\mu}_7$ and the associated completion $L^+_7$. The group $P_7$ is the hyperspecial parahoric subgroup stabilizing a certain Type 1 vertex of the building $X_7$ described in detail in \cite{CartwrightpadicC20}.

While we do not need such specific information, we briefly describe the lattice stabilized by $P_7$. One can find a unit $\eta \in L_7^\times$ such that $N_{L_7 / E_7}(\eta) = \de$ and $N_{L_7 / L^+_7}(\eta) = 1$. The element
\begin{equation}\label{eq:eta20}
c_{\eta} = \begin{pmatrix} \eta \phi(\eta) & 0 & 0 \\ 0 & \phi(\eta) & 0 \\ 0 & 0 & 1 \end{pmatrix}
\end{equation}
centralizes $\iota(x)$ for $x \in L$ with $\iota$ as in Equation~\eqref{eq:C20Emb1} and conjugates $\iota(y)$ to
\begin{equation*}%\label{eq:C20Emb3}
\begin{pmatrix} 0 & \eta & 0 \\ 0 & 0 & \phi(\eta) \\ \phi^2(\eta) & 0 & 0 \end{pmatrix}
\end{equation*}
with $\iota(y)$ as in Equation~\eqref{eq:C20Emb2}. Define
\begin{equation}\label{eq:Conjugator20}
c = \begin{pmatrix} 1 & 1 & 1 \\ \zeta - 1 & \phi(\zeta - 1) & \phi^2(\zeta - 1) \\ \al^{-1}(\zeta - 1)^2 & \al^{-1} \phi(\zeta - 1)^2 & \al^{-1} \phi^2(\zeta - 1)^2 \end{pmatrix}
\end{equation}
and the conjugate embedding
\begin{equation*}%\label{eq:New7Embed20}
\iota_\eta(z) = (c c_\eta) \iota(z) (c c_\eta)^{-1}
\end{equation*}
of $D$ into $\M_3(L_7)$. The hermitian element $h_c$ for this conjugate embedding is
\begin{align}
&\sig(c c_\eta)^{-1} h (c c_\eta)^{-1} \nonumber \\
=\,\, &\sig(c)^{-1} h c^{-1} \nonumber \\
=\,\, &\begin{pmatrix}
13 - 5 \al & 13 - 5 \al + 5 \beta - 2 \al \beta & \frac{-21 + 5 \al \beta}{2} + 4 \al - 7 \beta \\
13 - 5 \al - 5 \beta + 2 \al \beta & 16 - 6 \al & -14 + \frac{11 \al - 5 \beta}{2} + \al \beta \\
\frac{-21 - 5 \al \beta}{2} + 4 \al + 7 \beta & -14 + \frac{11 \al + 5 \beta}{2} - \al \beta & 14 - 5 \al
\end{pmatrix} \label{eq:20p7NewForm}
\end{align}
which we note is in fact a hermitian matrix in $\M_3(E)$. Then $P_7 < H_{F_7}$ is the open compact subgroup of consisting of the $g \in H_{F_7}$ such that $(c c_\eta) g (c c_\eta)^{-1}$ is $E_7$-adically integral. In other words, if $\calO_7$ is the ring of integral elements in $E_7$,
\begin{equation*}%\label{eq:Specific20P7}
P_7 \cong \GL_3(\calO_7) \cap (c c_\eta) H_{F_7} (c c_\eta)^{-1}.
\end{equation*}
It remains to define the open compact subgroups for the primes dividing $3$.

%%%%%%%%%%%%%%%%%%%%
\subsubsection{The open compact subgroups for $\frakp_{3\pm}$}\label{ssec:C20oc3}
%%%%%%%%%%%%%%%%%%%%

The construction uses several of the elements used to define $P_7$ in the previous subsection. We have $F_{3\pm} \cong \bbQ_3$ the completion associated with $\frakp_{3 \pm}$ and extensions $E_{3 \pm}$, $L_{3\pm}$ and $L^+_{3 \pm}$ all associated with the various extensions of $F$. When we only care about a particular field up to isomorphism, we drop the $\pm$ from the notation. We work in the product $H_{F_{3+}} \times H_{F_{3-}}$ and define maximal compact subgroups of each factor simultaneously.

A similar argument to $|\nu| = 7$, now where all fields are unramified extensions, shows that there exists an element $\eta \in L_3$ so that $N_{L_3 / E_3}(\eta) = \de$ and ${N_{L_3 / L^+_3}(\eta) = 1}$. From this $\eta$, we define $(c_\eta, c_\eta)$ analogously to Equation~\eqref{eq:eta20} and $c$ exactly as in Equation~\eqref{eq:Conjugator20}. We then define
\begin{align*}
c_+ &= \textrm{diag}(1,1,\pi_{3+}) \\ %\label{eq:c+20}
c_- &= \textrm{diag}(\pi_{3-},1,1) %\label{eq:c-20}
\end{align*}
and consider the two hermitian forms
\begin{equation*}%\label{eq:20p3NewForm}
h_{\pm} = \sig(c_\pm c c_\eta)^{-1} h (c_\pm c c_\eta)^{-1},
\end{equation*}
which, as with the form in Equation~\eqref{eq:20p7NewForm}, has matrix in $\M_3(E)$.

We first study $h_+$. Let $P_3^\prime$ be the elements in $H_{F_{3+}}$ whose conjugate by $c_+ c c_\eta $ preserve the standard $E_{3+}$ lattice and $P_3$ be the elements in $H_{F_{3-}}$ whose conjugate by $c_+ c c_\eta$ preserves the standard $E_{3-}$ lattice. Then $\pi_{3+} h_+$ has integral matrix entries over $E$ and determinant $2 \al - 5 \in \frakp_{3+}$, where $\nu_{3+}(2 \al - 5) = 1$ and $\nu_{3-}(2 \al - 5) = 0$. As described in \cite{CartwrightpadicC20}, this implies that the standard lattice in $E_{3+}$ is not self-dual for $h_+$ but its dual is a neighbor in the Bruhat--Tits building for $\PGL_3(E_3)$. However, $2 \al - 5$ is a unit of $F_{3-}$, which implies that the standard lattice is self-dual. These two cases precisely describe the two types of vertex in the Bruhat--Tits tree for $\SU(h_+, E_3)$, where $P_3^\prime$ stabilizes a Type 2 vertex and $P_3$ stabilizes a Type 1 vertex (see \S\ref{ssec:SUbuilding}).

For $h_-$, the roles are reversed. Indeed, it is now $\pi_{3-} h_-$ that is integral over $E$, and it has determinant $37 - 14 \al \in \frakp_{3-}$, which has $\nu_{3-}(37 - 14 \al) = 1$ and $\nu_{3+}(37 - 14 \al) = 0$. Recall from \S\ref{ssec:SUbuilding} that all Type $1$ (resp.~Type $2$) vertex stabilizers in $\SU(h_+, E_{3+}) \cong \SU(h_-, E_{3-})$ are conjugate, hence isomorphic. In particular, we can identify the elements in $H_{F_{3+}}$ whose conjugate by $c_- c c_\eta$ preserves the standard $E_{3+}$ lattice with $P_3$ and $P_3^\prime$ with the elements in $H_{F_{3-}}$ whose conjugate by $c_- c c_\eta$ preserves the standard $E_{3-}$ lattice.

%%%%%%%%%%%%%%%%%%%%
\subsubsection{The arithmetic lattices for $j = 47, 49$}
%%%%%%%%%%%%%%%%%%%%

We now define:
\begin{equation*}%\label{eq:47n49}
K_j = \begin{cases} \displaystyle{P_2 \times P_3^\prime \times P_3 \times \prod_{|\nu| \ge 5} P_\nu} & \quad j = 47 \\ \\ \displaystyle{P_2 \times P_3 \times P_3^\prime \times \prod_{|\nu| \ge 5} P_\nu} & \quad j = 49 \end{cases}
\end{equation*}
Here, the second factor of the product corresponds to the place of $F$ associated with $\frakp_{3+}$ and the third factor corresponds to the place for $\frakp_{3-}$. Then the image $\Gam_j$ of $\Gam(K_j)$ in $H^{ad}(F)^+$ is the fundamental group of the $j^{th}$ complex conjugate pair of fake projective planes. The groups $\Gam(K_{47})$ and $\Gam_{47}$ are the examples called $\Lam^\prime$ and $\conj{\Lam}{}^\prime$ in \cite[\S A.9]{PYFix}, and $\Gam(K_{49}), \Gam_{49}$ are $\Lam^{\prime \prime}, \conj{\Lam}{}^{\prime \prime}$ (this case was mistakenly omitted in \cite{PrasadYeung}, hence is not found there).

%%%%%%%%%%%%%%%%%%%%
\subsubsection{The isomorphism of profinite completions for $j = 47, 49$}
%%%%%%%%%%%%%%%%%%%%

%%%%%%%%%%%%%%%%%%%%
\begin{prop}\label{prop:47n49}
Let $X_{47}, X_{49}$ be fake projective planes with fundamental groups $\Gam_{47}, \Gam_{49}$ and choose a representative $\tau \in \Aut(\bbC)$ for the generator of $\Gal(F / \bbQ)$. Then, up to complex conjugation, $X_{47}^\tau$ is isomorphic to $X_{49}$. Therefore $\wh{\Gam}_{47} \cong \wh{\Gam}_{49}$.
\end{prop}
%%%%%%%%%%%%%%%%%%%%

%%%%%%%%%%%%%%%%%%%%
\begin{pf}
By Proposition~\ref{prop:MyForm}, the algebraic group associated with $X_{47}^\tau$ is the algebraic group with hermitian form $\tau(h)$, which is again $H$. By Theorem~\ref{thm:MilneSuh}, we must describe the open compact subgroup at each localization under the isomorphism of algebraic groups over the finite adeles of $F$. In particular, we must show that the self-isomorphism of $H(\bbA_F^\infty)$ induced by $\tau$ takes $K_{47}$ to $K_{49}$. By Proposition~\ref{prop:MyForm}, we only need to understand the $\tau$-conjugate of the local lattice used to defined each open compact subgroup $P_\nu$ of $H_{F_\nu}$.

For $|\nu| \neq 2,3,7$, $P_\nu$ is defined to be the stabilizer in $H_{F_\nu}$ of the standard $E_\mu$ lattice for $\mu$ dividing $\nu$. Since $\det(h)$ is a unit of $F$, the standard lattice is self-dual for $h$, $P_\nu$ is the stabilizer of a Type $1$ vertex of the building for $H_{F_\nu}$. Then $\tau(h)$ also has determinant a unit, so the stabilizer of the standard lattice is again $P_\nu$. Thus, up to conjugacy in $H_{F_{\tau(\nu)}}$ we send $P_\nu$ to $P_\nu \cong P_{\tau(\nu)}$ under the isomorphism. For $\frakp_7$, the determinant of the conjugate $h_c$ defined in Equation~\eqref{eq:20p7NewForm} is again a unit of $F$, so the same argument implies that $P_7$ is mapped to a conjugate of $P_7$. For $\frakp_2$, the open compact subgroup is defined by the unique maximal order of $D_{2+} \cong D_{2-}$, so it is again $P_2$ under the isomorphism induced by $\tau$.

It remains to consider $|\nu| = 3$. Here, we have that
\begin{equation*}%\label{eq:20v3swap}
\nu_{3 \pm}(\det(\tau(h))) = 1 - \nu_{3 \pm}(\det(h)).
\end{equation*}
This means that replacing $h$ with $\tau(h)$ swaps whether or not the standard lattice is self-dual with respect to the conjugate form described in the definition of the open compact subgroups for each of the two $3$-adic places in \S\ref{ssec:C20oc3}. It follows from Proposition~\ref{prop:MyForm} that $P_3^\prime \times P_3$ is mapped to $P_3 \times P_3^\prime$ by swapping factors under the Galois action of $\tau$. This combines with the assertions for the other $P_\nu$ to prove that the isomorphism of algebraic groups over $\bbA_F^\infty$ takes $K_{47}$ to $K_{49}$, which completes the proof of the proposition.
\end{pf}
%%%%%%%%%%%%%%%%%%%%

%%%%%%%%%%%%%%%%%%%%
\subsubsection{The isomorphism of profinite completions for $j = 48, 50$}
%%%%%%%%%%%%%%%%%%%%

%%%%%%%%%%%%%%%%%%%%
\begin{prop}\label{prop:48n50}
Let $X_{48}, X_{50}$ be fake projective planes with fundamental groups $\Gam_{48},\Gam_{50}$ and choose a representative $\tau \in \Aut(\bbC)$ for the generator of $\Gal(F / \bbQ)$. Then, up to complex conjugation, $X_{48}^\tau$ is isomorphic to $X_{50}$. Therefore $\wh{\Gam}_{48} \cong \wh{\Gam}_{50}$.
\end{prop}
%%%%%%%%%%%%%%%%%%%%

%%%%%%%%%%%%%%%%%%%%
\begin{pf}
For $j = 47, 49$, using computer algebra software, we see that there is a fake projective plane $X_j$ with fundamental group $\Gam_j$ such that there is a unique connected Galois \'etale cover $Y_j \to X_j$ with group $S_3$ so that $Y_j$ is also a Galois \'etale cover of $X_{j+1}$ with group $\bbZ / 6$. By Lemma~\ref{lem:AutProperties}, the map $Y_{47} \to Y_{47}^\tau$ induces $X_{47} \to X_{47}^\tau$. Proposition~\ref{prop:47n49} implies that, up to replacing $X_{49}$ with its complex conjugate, $X_{47}^\tau \cong X_{49}$, so $Y_{47}^\tau \cong Y_{49}$ by uniqueness of each $Y_j$ as a $6$-fold Galois cover of $X_j$ that also covers $X_{j+1}$. We then obtain the diagram of maps in Figure~\ref{fig:C20etale}, where Corollary~\ref{cor:FakeQuotient} implies the existence of the dashed arrow. Therefore, up to the choice of complex structure on $X_{50}$ determined by the choice of complex structure on $X_{49}$, this proves that $X_{50} \cong X_{48}^\tau$, which proves the proposition.
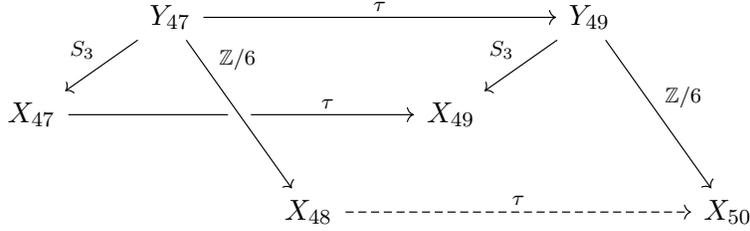
\begin{figure}[h]
\centering
\begin{tikzcd}
& Y_{47} \arrow[rrr, "\tau"] & & & Y_{49} \arrow[ddr, "\bbZ / 6"] & \\
X_{47} \arrow[<-, ur, "S_3"] \arrow[rrr, "\tau" near end] & & & X_{49} \arrow[<-, ur, "S_3"] & \\
& & X_{48} \arrow[rrr, dashed, "\tau"] \arrow[<-, uul, crossing over, "\bbZ / 6" above right, near end] & & & X_{50}
\end{tikzcd}
\caption{The cover $Y_j$ for $j = 47, 49$.}\label{fig:C20etale}
\end{figure}
\end{pf}
%%%%%%%%%%%%%%%%%%%%

%%%%%%%%%%%%%%%%%%%%
\section*{Appendix: Fake projective plane groups}
%%%%%%%%%%%%%%%%%%%%

We label the fifty (topological) fundamental groups of fake projective planes following notation used by Cartwright--Steger \cite{CartwrightSteger} and Prasad--Yeung \cite{PrasadYeung, PYFix}. First, fake projective planes are grouped by the totally imaginary field used to define the relevant commensurability classes of arithmetic lattices. Within each commensurability class, we number the lattices based on their order of appearance in the register of presentations of fake projective plane groups given on Cartwright's website \cite{Cartwright}; this is the file that gives each fake projective plane group the `identifier' given below. Specifically, each table below consists of a number $j$ and the name from \cite{Cartwright}, meaning that throughout this paper $\Gam_j$ denotes the fundamental group of the given complex conjugate pair of fake projective planes.

\begin{table}[h!]
\centering
\begin{tabular}{|c|l|}
\hline
$j$ & \multicolumn{1}{|c|}{Identifier from \cite{Cartwright}} \\
\hline
$1$ & $\, p=5,\emptyset,D_3$ \\
$2$ & $\, p=5,\{2\},D_3$ \\
$3$ & $\, p=5,\{2i\}$ \\
\hline
\end{tabular}
\caption{The cases $a = 1$ from \cite{Cartwright}: $\bbQ(\sqrt{-1})$}
\end{table}

\begin{table}[h!]
\centering
\begin{tabular}{|c|l|}
\hline
$j$ & \multicolumn{1}{|c|}{Identifier from \cite{Cartwright}} \\
\hline
$4$ & $\, p=3,\emptyset,D_3$ \\
$5$ & $\, p=3,\{2\},D_3$ \\
$6$ & $\, p=3,\{2i\}$ \\
\hline
\end{tabular}
\caption{The cases $a = 2$ from \cite{Cartwright}: $\bbQ(\sqrt{-2})$}
\end{table}

\begin{table}[h!]
\centering
\begin{tabular}{|c|l|c|l|}
\hline
$j$ & \multicolumn{1}{|c|}{Identifier from \cite{Cartwright}} & $j$ & \multicolumn{1}{|c|}{Identifier from \cite{Cartwright}} \\
\hline
$7$ & $\, p=2,\emptyset,D_3\, 2_7$ & $14$ & $\, p=2,\{3\},D_3$ \\
$8$ & $\, p=2,\emptyset,7_{21}$ & $15$ & $\, p=2,\{3\},3_3$ \\
$9$ & $\, p=2,\emptyset,D_3\, X_7$ & $16$ & $\, p=2,\{3,7\},D_3$ \\
$10$ & $\, p=2,\{7\},D_3\, 2_7$ & $17$ & $\, p=2,\{3,7\},3_3$ \\
$11$ & $\, p=2,\{7\},D_3\, 7_7$ & $18$ & $\, p=2,\{5\}$ \\
$12$ & $\, p=2,\{7\},D_3 7^\prime_7$ & $19$ & $\, p=2,\{5,7\}$ \\
\cline{3-4}
$13$ & $\, p=2,\{7\},7_{21}$ & \multicolumn{1}{c}{} & \multicolumn{1}{c}{} \\
\cline{1-2}
\end{tabular}
\caption{The cases $a = 7$ from \cite{Cartwright}: $\bbQ(\sqrt{-7})$}
\end{table}

\begin{table}[h!]
\centering
\begin{tabular}{|c|l|c|l|}
\hline
$j$ & \multicolumn{1}{|c|}{Identifier from \cite{Cartwright}} & $j$ & \multicolumn{1}{|c|}{Identifier from \cite{Cartwright}} \\
\hline
$20$ & $\, p=2,\emptyset,D_3$ & $25$ & $\, p=2,\{5\},D_3$ \\
$21$ & $\, p=2,\emptyset,3_3$ & $26$ & $\, p=2,\{5\},3_3$ \\
$22$ & $\, p=2,\{3\},D_3$ & $27$ & $\, p=2,\{3,5\},D_3$ \\
$23$ & $\, p=2,\{3\},3_3$ & $28$ & $\, p=2,\{3,5\},3_3$ \\
$24$ & $\, p=2,\{3\},(D3)_3$ & $29$ & $\, p=2,\{3,5\},(D3)_3$ \\
\hline
\end{tabular}
\caption{The cases $a = 15$ from \cite{Cartwright}: $\bbQ(\sqrt{-15})$}
\end{table}

\begin{table}[h!]
\centering
\begin{tabular}{|c|l|}
\hline
$j$ & \multicolumn{1}{|c|}{Identifier from \cite{Cartwright}} \\
\hline
$30$ & $\quad p=2,\emptyset$ \\
$31$ & $\quad p=2,\{23\}$ \\
\hline
\end{tabular}
\caption{The cases $a = 23$ from \cite{Cartwright}: $\bbQ(\sqrt{-23})$}
\end{table}

\begin{table}[h!]
\centering
\begin{tabular}{|c|l|c|l|}
\hline
$j$ & \multicolumn{1}{|c|}{Identifier from \cite{Cartwright}} & $j$ & \multicolumn{1}{|c|}{Identifier from \cite{Cartwright}} \\
\hline
$32$ & $\quad p=2,\emptyset,d_3\, D_3$ & $36$ & $\quad p=2,\emptyset,d_3\, X^\prime_3$ \\
$33$ & $\quad p=2,\emptyset,D_3\, X_3$ & $37$ & $\quad p=2,\emptyset,X_9$\\
$34$ & $\quad p=2,\emptyset,(dD)_3\, X_3$ & $38$ & $\quad p=2,\{3\},d_3\, D_3$ \\
\cline{3-4}
$35$ & $\quad p=2,\emptyset,(d^2D)_3\, X_3$ & \multicolumn{1}{c}{} & \multicolumn{1}{c}{} \\
\cline{1-2}
\end{tabular}
\caption{The cases $\calC_2$ from \cite{Cartwright}: $\bbQ(\sqrt{5}, e^{2 \pi i / 3})$}\label{tb:C2}
\end{table}

\begin{table}[h!]
\centering
\begin{tabular}{|c|l|}
\hline
$j$ & \multicolumn{1}{|c|}{Identifier from \cite{Cartwright}} \\
\hline
$39$ & $\quad p=2,\emptyset,D_3$ \\
$40$ & $\quad p=2,\{17-\},D_3$ \\
\hline
\end{tabular}
\caption{The cases $\calC_{10}$ from \cite{Cartwright}: $\bbQ(\sqrt{-7+4 \sqrt{2}})$}
\end{table}

\begin{table}[h!]
\centering
\begin{tabular}{|c|l|}
\hline
$j$ & \multicolumn{1}{|c|}{Identifier from \cite{Cartwright}} \\
\hline
$41$ & $\quad p=3,\emptyset,d_3\, D_3$ \\
$42$ & $\quad p=3,\{2\},D_3$ \\
$43$ & $\quad p=3,\{2\},(dD)_3$ \\
$44$ & $\quad p=3,\{2\},(d^2D)_3$ \\
$45$ & $\quad p=3,\{2i\}$ \\
\hline
\end{tabular}
\caption{The cases $\calC_{18}$ from \cite{Cartwright}: $\bbQ(\sqrt{6}, e^{2 \pi i / 3})$}\label{tb:C18}
\end{table}

\begin{table}[h!]
\centering
\begin{tabular}{|c|l|}
\hline
$j$ & \multicolumn{1}{|c|}{Identifier from \cite{Cartwright}} \\
\hline
$46$ & $\quad p=2,\emptyset,D_3\, 2_7$ \\
$47$ & $\quad p=2,\{3+\},D_3$ \\
$48$ & $\quad p=2,\{3+\},(3+)_3$ \\
$49$ & $\quad p=2,\{3-\},D_3$ \\
$50$ & $\quad p=2,\{3-\},(3-)_3$ \\
\hline
\end{tabular}
\caption{The cases $\calC_{20}$ from \cite{Cartwright}: $\bbQ(\sqrt{7}, i)$}
\end{table}

\pagebreak

%%%%%%%%%%%%%%%%%%%%
\bibliography{FakeP2pi1}
%%%%%%%%%%%%%%%%%%%%

%%%%%%%%%%%%%%%%%%%%
\end{document}